\documentclass[12pt]{article}
\input epsf
\usepackage{epsfig}
\usepackage{latexsym}
\usepackage{amsfonts}
\usepackage{amsmath}
\usepackage[psamsfonts]{amssymb}
\usepackage{float}
\def\R{\mathbf{R}}
\def\C{\mathbf{C}}
\def\Z{\mathbf{Z}}
\def\Tr{\mathrm{Tr}\,}

\def\mod{\mathrm{mod}\,}
\def\Aut{\mathrm{Aut}\,}
\title{Circular pentagons and real solutions of Painlev\'e VI equations}
\author{Alexandre Eremenko\thanks{Supported by NSF grant DMS--1361836.}\; and
Andrei Gabrielov\thanks{Supported by NSF grant DMS--1161629.}
}
\begin{document}
\maketitle
\begin{abstract}
We study real solutions of a class of Painlev\'e VI equations.
To each such solution we associate a geometric object,
a one-parametric family of circular pentagons.
We describe an algorithm which permits to compute the numbers of zeros,
poles, $1$-points and fixed
points of the solution on the interval $(1,+\infty)$ and their mutual position.
The monodromy of the associated
linear equation and parameters of the Painlev\'e VI equation
are easily recovered from the family of pentagons.
\vspace{.1in}

\end{abstract}

\noindent
{\bf 0. Introduction}
\vspace{.1in}

Consider a linear differential equation
\begin{equation}\label{1}
w''-Pw'+Qw=0,
\end{equation}
where $P$ and $Q$ are rational functions of the complex
independent variable,
and assume that all singularities
are regular, and all parameters (singular points,
exponents and accessory
parameters) are real.
Then the ratio $f=w_1/w_2$ of two linearly independent
solutions maps the upper half-plane
onto a circular polygon (see Section 3 below for a precise definition).
Every simply connected circular polygon can be obtained this way.
Klein \cite{Klein} and Van Vleck \cite{VV2}
used this connection between differential equations with three
singularities and circular triangles to count zeros and poles of hypergeometric
functions. We use a similar idea to count special points on
a real interval of real solutions
of Painlev\'e VI equations
\begin{eqnarray}\label{PVI}
q_{xx}&=&\frac{1}{2}\left(\frac{1}{q}+\frac{1}{q-1}+\frac{1}{q-x}\right)q_x^2
-\left(\frac{1}{x}+\frac{1}{x-1}+\frac{1}{q-x}\right)q_x \\ \nonumber
&+&\frac{q(q-1)(q-x)}{2x^2(x-1)^2}\left\{\kappa_4^2-\kappa_1^2\frac{x}{q^2}
+\kappa_2^2\frac{x-1}{(q-1)^2}+(1-\kappa_3^2)\frac{x(x-1)}{(q-x)^2}\right\},
\end{eqnarray}
where $q$ is a function of $x$, and $\kappa_j$ are real parameters.
By definition, {\em special points} of a solution are those
points $x$ for which $q(x)\in\{0,1,x,\infty\}$,
the points where the
assumptions of the
existence and uniqueness theorem of Cauchy are violated.

Equation (\ref{PVI}), which we call PVI, was originally discovered by
Picard \cite{Pic}, Painlev\'e \cite{Pai} and
Gambier \cite{Gam} as the most
general equation of the form
\begin{equation}\label{2de}
q_{xx}=R(q_x,q,x),
\end{equation}
where $R$ is a rational function
of $q_x,q$ whose coefficients are analytic in $x$,
and whose solutions have no movable singularities
(poles are not counted as singularities).
For PVI, this means that all solutions admit a meromorphic continuation
in the region $\C\backslash\{0,1\}$. Painlev\'e and Gambier
classified all equations (\ref{2de}) without movable singularities,
and found
that all except six of them can be reduced to linear or first order
differential equations.
Of the six remaining equations,
PVI is the most general, in the sense
that the other
five
can be obtained from it by certain degeneration process.

Meanwhile, Richard Fuchs
\cite{Fuchs}
independently discovered (\ref{PVI}) as the condition of isomonodromic
deformation of a linear differential equation (\ref{1})
with five regular singularities,
one of them apparent with exponents $0$ and $2$.
Later it was found that Painlev\'e equations arise
in a variety of problems of mathematics and physics
\cite{Dubrovin,Hitchin,KM,Z}, and their solutions,
called Painlev\'e transcendents,
gradually gain the status of special functions.

Besides applications, most work on PVI falls into three categories:
a) transformations of the equation \cite{Okamoto,inaba}, b) search
for special solutions, like algebraic ones \cite{Lisovyy} or those
expressed in terms of classical special functions \cite{Hitchin}, and
c) asymptotics at the fixed singularities $0,1,\infty$, \cite{Jimbo,Guzzetti}.

In this paper, we study real solutions of PVI with real
parameters $\kappa_j$
on a real interval, one of the three intervals between the
fixed singularities $(0,1,\infty)$. Our main result is a combinatorial
algorithm
(which can be performed without a computer)
which determines the number of special
points and their mutual position on the interval.
The outcome of the algorithm is a sequence of symbols
$0,1,x,\infty$ which shows in which order the solution
$q(x)$ of (\ref{PVI}) takes these four values as $x$ increases
from $1$ to $\infty$. In particular, we describe those solutions
which have no special points on $(1,\infty)$.

So we obtain
global, exact (non-asymptotic) results describing qualitative features
of a quite general class of solutions, namely real solutions.

We do this by exploring the connection with a linear differential
equation of the type (\ref{1}) with 5 regular singularities discovered by
Fuchs.
When all parameters of this linear equation
are real, the ratio $f$ of any two linearly independent solutions
of (\ref{1}) maps the upper
half-plane onto a circular pentagon of a special kind: it is a circular
quadrilateral with a slit or a circular triangle with a slit.
If the $f$-preimage of the tip of the slit is not counted as a vertex,
our special pentagon can be considered as a conformal quadrilateral.
A real solution of PVI describes the relation
between the modulus of this conformal quadrilateral
and the $f$-preimage of the tip of the slit.
For some values of the modulus the slit vanishes or becomes
a cross cut, and the pentagon becomes a circular quadrilateral.
These values of
the modulus correspond to special points of the PVI solution.
Thus the study of the number and mutual position of these special points
is equivalent to a geometric problem of describing the evolution of
a one-parameter family of special pentagons.
We do this using a representation of
circular polygons by cell decompositions of a disk
which are called nets. This method was developed in
\cite{EGA,EGSV,EG2,EGT1,EGT2}
for other problems.
\vspace{.1in}

We are grateful to C.-S. Lin who shared with us the result
from his unpublished preprint \cite{chen2} which stimulated
this paper. We also thank Philip Boalch, Alexander Its, Oleg Lisovyy,
and Vitaly Tarasov
for illuminating discussions of PVI.
\vspace{.05in}

\noindent
{\bf 1. A class of linear differential equations}
\vspace{.1in}

We consider the class of linear differential equations
(\ref{1}) of the form
\begin{equation}\label{1a}
w''-\left(\frac{1}{z-q}+\sum_{j=1}^3\frac{\kappa_j-1}{z-t_j}\right)w'
+\left(\frac{p}{z-q}-\sum_{j=1}^3\frac{h_j}{z-t_j}\right)w=0,
\end{equation}
where
\begin{equation}\label{ts}
(t_1,t_2,t_3)=(0,1,x),
\end{equation}
$x>1,\;\kappa_j\geq 0,\; p$ and $q$ are real numbers.
We impose the following conditions:

a) $\infty$ is a regular singularity, with exponent difference
$\kappa_4\geq 0$,

b) $q$ is an apparent singularity,
but the singularities at $(0,1,x,\infty)$
are not apparent (have non-trivial local projective monodromy).

It follows from the form of (\ref{1a}) that the exponents
at $q$ are $0$ and $2$, and $t_j$ are regular singularities
with exponents $0$ and $\kappa_j$
for $1\leq j\leq 3$. The exponents at $\infty$ are determined
from the Fuchs relation and from the assumption that their difference
is $\kappa_4$.

Conditions a), b) determine the $h_j$ uniquely in terms of
$p,q,x$ and $\kappa_j, 1\leq j\leq 4$,
by solving the
following system of linear equations:
\begin{eqnarray}\label{a}
\sum_{j=1}^3h_j&=&p,\\ \label{b}
\sum_{j=1}^3t_jh_j&=&pq+\frac{\kappa_4^2-1}{4}-
\frac{1}{4}\left(\sum_{j=1}^3\kappa_j\right)^2+
\frac{1}{2}\sum_{j=1}^3\kappa_j,\\ \label{c}
\sum_{j=1}^3\frac{h_j}{q-t_j}&=&p^2-p\sum_{j=1}^3\frac{\kappa_j-1}{q-t_j}.
\end{eqnarray}
Equations (\ref{a}) and (\ref{b}) correspond to condition a),
while equation (\ref{c}) corresponds to condition b).

The determinant of this system is
$$\frac{x(x-1)}{q(q-1)(q-x)},$$
thus for given real $\kappa_j\geq 0,\;p,\; x\not\in\{0,1\}$,
and $q\not\in\{0,1,x\}$
the coefficients $h_j$ are uniquely defined real numbers.
Our notation in
(\ref{PVI}) and (\ref{1a}) is the same as in \cite{inaba},
but we notice a misprint in the first line of \cite[(2.1)]{inaba}:
$q$ must be~$q_x$.
\vspace{.1in}

\noindent
{\bf 2. Isomonodromic deformation and PVI}
\vspace{.1in}

\begin{figure}
\centering
\includegraphics[width=2.5in]{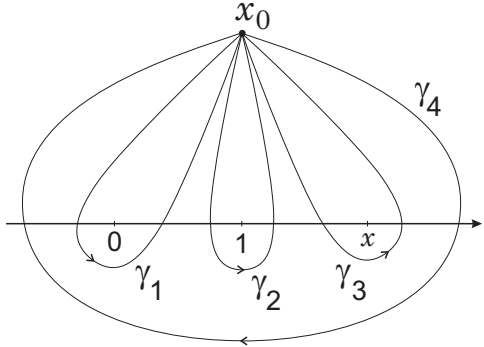}
\caption{Loops defining $M_1,\,M_2,\,M_3$, and $M_4$.}\label{loops}
\end{figure}

Let us choose the
generators of the fundamental
group of $\Omega=\C\backslash\{0,1,x\}$
as shown in Fig.~\ref{loops}.

Let $(w_1,w_2)$ be a pair
of linearly independent solutions
of (\ref{1a}) normalized by
\begin{equation}\label{normalization}
\left(\begin{array}{cc}w_1(x_0)& w_2(x_0)\\
w_1^\prime(x_0)& w_2^\prime(x_0)\end{array}\right)=I.
\end{equation}
Performing an analytic continuation of these solutions
along an element $\gamma\in\pi_1(\Omega,x_0)$ we obtain
$$(w_1^\gamma,w_2^\gamma)=(w_1,w_2)M_\gamma$$
for some $M_\gamma\in GL(2,\C)$.
Notice that the map $\gamma\mapsto M_\gamma$ is an anti-representation
of the fundamental group\footnote{For the fundamental group we use the
standard notation: product $\gamma_1\gamma_2$ means that the path
$\gamma_1$ is followed by the path $\gamma_2$.}

For the ratio $f=w_1/w_2$
we obtain
$$f^\gamma=T_\gamma\circ f,$$
where $T_\gamma$ is a linear-fractional transformation.
We identify the group of linear-fractional transformations
with $PSL(2,\C)=SL(2,\C)/\{\pm I\}$,
and the quadruple $(T_1, T_2, T_3, T_4)$ is the set of
generators of the {\em projective monodromy representation}
$\pi_1(\Omega,x_0)\to PSL(2,\C)$. The correspondence $\gamma\to T_\gamma$
is a group homomorphism.
The generators
are chosen so that
\def\id{\mathrm{id}}
\begin{equation}\label{id}
T_1 T_2 T_3 T_4 =\id,
\end{equation}
and we assume that
\begin{equation}\label{as}
T_j\neq \id,\quad 1\leq j\leq 4.
\end{equation}
For the matrices representing the generators $T_j$ we use the same letters,
and they are related to the matrices $M_j$ by
$$T_j=\frac{1}{\sqrt{\det M_j}}M_j^T,\quad 1\leq j\leq 4,$$
where $T$ stands for the transposition.

When the parameters $\kappa_j$ are fixed, the projective
monodromy representation of equation (\ref{1})
depends on $x$, $p$ and $q$.

When we change $x$ and deform the loops continuously,
the condition that the monodromy matrices do not change is that
$p=p(x)$ and $q=q(x)$ satisfy the following non-autonomous
Hamiltonian system
\cite[(3.7)]{inaba}:
$$\frac{dq}{dx}=\frac{\partial h}{\partial p},\quad
\frac{dp}{dx}=-\frac{\partial h}{\partial q}.$$
Here the Hamiltonian $h=h_3$ is given by\footnote{See (\ref{1a})
for the definition of $h_3$.}
\begin{eqnarray*}
x(x-1)h&=&q(q-1)(q-x)p^2\\
&-&\left\{(\kappa_3-1)q(q-1)+\kappa_1(q-1)(q-x)
+
\kappa_2q(q-x)\right\}p\\
&+&\kappa_0(\kappa_0+\kappa_4)(q-x),
\end{eqnarray*}
where
$$\kappa_0=(1-\kappa_1-\kappa_2-\kappa_3-\kappa_4)/2.$$
This Hamiltonian system is equivalent to (\ref{PVI}).
All solutions of PVI are obtained in this way.

Special points of $q(x)$ correspond to collisions of
the singular point $q$ with one of the four other
singular points of equation (\ref{1a}).
Thus when $x$ is a special point,
(\ref{1}) becomes an equation with four regular singularities
(Heun's equation).


In this paper we consider parameters $\kappa_j\geq 0,\; 1\leq j\leq 4,$
and real solutions $q(x)$ of (\ref{PVI}) defined for $x\in(1,+\infty)$.
In view of the formulas (\ref{a}), (\ref{b}), (\ref{c}),
in this case all parameters in (\ref{1a}) are real.

The condition on the monodromy matrices that ensures that the solution of
PVI is real is discussed in Appendix I.
\vspace{.1in}

\noindent
{\em Remark.} A more general class of real equations (\ref{PVI})
is obtained by allowing some $\kappa_j$ be pure imaginary.
In this case, equation (\ref{1a}) also has a geometric interpretation
\cite{Yo1}, but very different from the interpretation
in this paper: the developing map $f$ (defined in the next section)
still maps the upper half-plane onto a Riemann
surface bounded by four circles,
but when some $\kappa_j$ are imaginary, this Riemann surface has infinitely
many sheets. 
\vspace{.1in}

\noindent
{\bf 3. Circular polygons}
\vspace{.1in}

{\em A circular $n$-gon}
is a bordered surface homeomorphic to a closed disk,
spread over the sphere without
ramification points in the interior, and such that the border consists
of $n$ arcs and $n$ points separating them, so that
each arc projects into a circle on the sphere locally injectively.

To give a more formal definition, we denote by $S$ the conformal sphere
(the unique compact simply connected Riemann surface).
A {\em circle} in $S$
can be defined by using only the conformal structure:
it is the set of fixed
points of an anti-conformal involution. Conformal automorphisms of $S$
send circles to circles.

Let $\overline{D}$ be a conformal closed disk\footnote{The closure of a Jordan region.}, and let $\{ t_j\}$ be $n$ distinct
boundary points enumerated according to the standard orientation
of $\partial D$. In what follows we understand the subscript
$j$ in $t_j$ and in other similar notations as a residue modulo $n$.

A {\em developing map} is a
continuous function $f:\overline{D}\to S$ which is
holomorphic in $\overline{D}\backslash\{ t_1\ldots,t_n\}$ and satisfies
\begin{equation}\label{crit}
f'(z)\neq 0,\quad z\in\overline{D}\backslash\{ t_1,\ldots,t_n\},
\end{equation}
\begin{equation}\label{loc1}
f(z)\sim f(t_j)+c_j(z-t_j)^{\alpha_j}\quad\mbox{as}\quad z\to t_j,
\end{equation}
where $\alpha_j>0$ and $c_j\in\C^*$, or
\begin{equation}\label{log}
f(z)\sim f(t_j)+c_j/\log(z-t_j)\quad\mbox{as}\quad z\to t_j,
\end{equation}
and such that $f([t_{j-1},t_j])$ are contained in some circles $C_j\subset S$.
(These formulas need an evident modification if
$f(t_j)=\infty$, or $f(z)=\infty$ in (\ref{crit}).)
The circles $C_j$ need not be distinct.
A {\em circular $n$-gon} is identified with the ordered set
\begin{equation}\label{ngon}
(\overline{D},t_1,\ldots,t_n,f).
\end{equation}
Sometimes we will omit the word ``circular'',
calling these objects simply polygons (digons, triangles, quadrilaterals, etc.)

Two circular polygons
\begin{equation}\label{two}
(\overline{D}_1,t_1^\prime,\ldots,t_n^\prime,f_1)\quad\mbox{and}\quad
(\overline{D}_2,t_1^{\prime\prime},\ldots,t_n^{\prime\prime},f_2)
\end{equation}
are considered {\em equal}
if there exists a conformal map
$\phi:\overline{D}_1\to \overline{D}_2$ such
that $\phi(t_j^\prime)=t_j^{\prime\prime}$
and $f_1=f_2\circ\phi.$ If the last equality is replaced
by $f_1=L\circ f_2\circ\phi$, where $L\in\Aut{S}$ then the two polygons are
called {\em equivalent}.

The points $t_j$ are called {\em corners} and the arcs $(t_{j-1},t_j)$
{\em sides} of a polygon. The angle at $t_j$ is defined as $\alpha_j$ in (\ref{loc1}),
and we set $\alpha_j=0$ if (\ref{log}) holds.

{\em Notice that we measure the
angles in half-turns rather than radians.}

We denote by $C_j$ the circle containing $f([t_{j-1},t_j])$. Then we
obtain $n$ labeled circles with the property that
\begin{equation}\label{C1}
C_{j}\cap C_{j+1}\neq\emptyset,\quad j\in\Z_4.
\end{equation}
Indeed, $f(t_j)$ belongs to this intersection.
Any such sequence of circles will be called
an $n$-{\em circle chain}, or simply a chain when it is clear what $n$ is.

Notice that $0$-gons are just disks,
while $1$-gons are disks with
one marked point where the angle is $1$.

Sometimes it will be convenient to use a Riemannian
metric on our polygons.
To introduce it, start with the standard spherical Riemannian
metric of curvature $1$
on $S$ and pull it back to $\overline{D}$ via $f$.
The resulting metric $\rho$ on $\overline{D}$ is a conformal Riemannian metric
of curvature $1$ on $\overline{D}\backslash\{ t_1,\ldots,t_n\}$, has
conic singularities with the angles $\alpha_j$ at $t_j$,
and each side $(t_j,t_{j+1})$ has constant geodesic curvature.
All metric spaces with these properties
arise from circular polygons.

In what follows the word ``distance''
will always mean {\em intrinsic distance}:
the infimum of lengths of curves connecting two points,
where the length of a curve is measured using the intrinsic
metric.
The area of an $n$-gon is also measured in the
pull-back spherical metric. It is easy to see that all our polygons
have finite areas, moreover, the preimages of points under
developing maps are finite.

Polygons are equal if and only if the corresponding
metric spaces are isometric.
Of course, equivalent polygons may be different as metric
spaces\footnote{One could use only $PSL(2,\C)$-invariant
notions, like cross-ratios instead of distances etc., as Klein does.
But we find the metric notions more convenient and more intuitive.}.
\vspace{.01in}

{\em Gluing of two polygons.}
\vspace{.1in}

We will use the operation of {\em gluing} circular polygons along a
``matching'' boundary arc. Suppose that for two polygons in (\ref{two}),
$D_1$ and $D_2$ are the upper and lower halves of the unit disk,
the interval $(-1,1)$ contains no corners of either polygon, and
is mapped by $f_1$ and $f_2$ to the same arc of a circle\footnote{This means
that there is an increasing diffeomorphism $\psi:[-1,1]\to[-1,1]$ such that
$f_2(t)=f_1\circ\psi(t),\; t\in[-1,1]$.}
(This ``arc''
can be longer than the whole circle.)
Then there exist a simple
curve $\gamma$ in the unit disk $D$ with endpoints at $\pm1$ dividing
$D$ into two regions $D'$ and $D''$ and conformal
homeomorphisms $\phi_1:D'\to D_1$ and $\phi_2:D''\to D_2$ such that
$\phi_j(\pm1)=\pm1$ and $f_1\circ\phi_1(z)=f_2\circ\phi_2(z),\; z\in[0,1]$.
Such conformal homeomorphisms $\phi_j$ exist by a theorem of Lavrentiev,
\cite[Ch. VI, \S1]{Goluzin}.
Then
$$f(z)=\left\{\begin{array}{ll}f_1\circ\phi_1(z),& z\in D'\\
f_2\circ\phi_2(z),& z\in D''\end{array}\right.$$
extended by continuity on $\gamma$,
is the developing map of a new polygon which is called the gluing of
our two polygons along the common boundary arc.
\vspace{.1in}

{\em Variation of the slit.}
\vspace{.1in}

Consider an $(n+1)$-gon $Q=(\overline{D},t_1,\ldots,t_n;q,f)$,
where the corner $q$ can be anywhere between the $t_j$,
this is why we use a different notation for this corner.
Suppose that the angle at $q$ equals $2$ and the $f$-images of
the two sides
meeting at $q$ belong to the same circle $C$. (If $q\in(t_{k-1},t_k)$
then $C=C_k$.)
This means that $f$ maps a small neighborhood $V$
of $q$ in $D$ homeomorphically onto a disk centered at $f(q)$ with a slit
from the center to the circumference along an arc of the circle
$C$.

In this situation
we say that the polygon
has a slit, and
$b:=f(q)$ is called the {\em tip of the slit}.
The slit itself is formally defined as follows:

{\em The slit is the maximal interval $[t,t']$ such that
$q\in[t,t']\subset[t_{k-1},t_k]$ and
the intrinsic lengths of $[t,q]$ and $[q,t']$ are equal.}

Consider the small arc
$\gamma\in D$ with endpoints $q$ and $c\in D$ which is defined
by $\gamma=f^{-1}(C)\cap V$.

Let $\phi$ be a conformal map
of $D$ onto $D\backslash\gamma$. Then $f_1=f\circ\phi$
defines a new $(n+1)$-gon with
corners $t^\prime_j=\phi^{-1}(t_j)$ for $1\leq j\leq n$
and
$q':=\phi^{-1}(c)$. We say that this new polygon is obtained from
the old one by {\em lengthening the slit}, and the old polygon is obtained
from the new one by {\em shortening the slit}.

Here is an alternative explanation of variation of the slit.
Suppose that $D=H$ and $q\in(t_{k-1},t_k)$. As the sides $(t_{k-1},q)$
and $(q,t_k)$ are mapped by $f$ to the same circle $C$, we can extend
$f$ by reflection to the lower half-plane $H^*$. The resulting function
$\tilde{f}$ is meromorphic in $G=H\cup H^*\cup(t_{k-1},t_k)$, maps
$(t_{k-1},t_k)$ into a circle $C$ and has
exactly one simple critical point at $q$. Let $\sigma$ be the reflection
in $C$. Choose a small disk $B$ centered at $f(q)$ and let $V\subset G$
be the component of $\tilde{f}^{-1}(B)$ which contains $q$.
Let $\psi_s$
be a family of diffeomorphisms of the sphere $S$, which commutes with $\sigma$,
whose restriction to $S\backslash B$ is the identity map, and which moves
$\tilde{f}(q)$ to a point $s\in C$ near $f(q)$. Then the main existence
theorem for quasiconformal mappings in \cite{Ahl2}
implies that there is a quasiconformal homeomorphism $\phi_s:G\to G$
which commutes with
complex conjugation and such
that $f_s=\psi_s\circ \tilde{f}\circ\phi_s$ is holomorphic.
The restriction of $f_s$ onto $H$ is the
developing map of the deformed polygon
$Q'=(\overline{H},t^\prime_1,\ldots,t^\prime_4; q',
f_s)$ where $t_j^\prime=\phi^{-1}_s(t_j)$ and $q'=\phi^{-1}_s(q)$.
The dependence of $f_s$ on $s$ is real analytic.
\vspace{.1in}

Whenever we have a slit it can be lengthened or shortened.
This operation does not affect the angles, the chain of the polygon,
or the images of the sides other than those two meeting at $q$.
\vspace{.1in}

\noindent
{\bf 4. Relation between equation (\ref{1a}) and a class of circular
pentagons}
\vspace{.02in}
\nopagebreak

We consider equation (\ref{1a})
satisfying conditions a) and b) after (\ref{1a}).
\vspace{.1in}

\noindent
{\bf Proposition 1.} {\em If $p\in\R$, $x>1$, and
$q\in\R\backslash\{0,1,x\}$, then the ratio $f$
of any two linearly
independent solutions of (\ref{1a})
is the developing map of a circular
pentagon with $D=H$,
and corners at $(t_1,t_2,t_3,t_4)=(0,1,x,\infty)$ and $q$,
with the angles $\kappa_j$ at $t_j$, $1\leq j\leq 4$,
and $2$ at $q$. The $f$-images of the two sides meeting
at $q$ belong to the same circle.

Conversely, the developing map of every circular pentagon with
such properties
is the ratio of two solutions of an equation (\ref{1a}) satisfying
conditions a), b),
with all parameters real, and $0,1,x,q$ all distinct.

The projective monodromy group of (\ref{1a})
consists of the products of
even numbers of reflections in the sides of the pentagon.
Condition (\ref{as}) holds if and only if no pair of sides meeting
at $t_j,\; 1\leq j\leq 4,$
is mapped by the developing map into the same circle.}
\vspace{.1in}

Notice that $q$ can be on any of the four intervals
$(t_{j-1},t_j),\; j\in\Z_4$.
\vspace{.1in}

{\em Proof.} Let $f=w_1/w_2$ be the ratio of linearly independent solutions.
Then $f'=(w_1^\prime w_2-w_1w_2^\prime)/w_2^2$, so $f$ is locally univalent
in the upper half-plane. If we impose real initial conditions at some real
non-singular point, both solutions will be real, and $f$
will be real on the interval between the singularities containing this point.
Any other initial condition will give new $f$ related to the old one
by a linear-fractional transformation, so $f(z)$ maps
every interval between the singular points onto an arc of a circle.
The exponents at a singular
point $t_j$ are $0$ and $\kappa_j$ if $1\leq j\leq 3$, so locally $f(z)$
behaves as in (\ref{loc1}), (\ref{log}).
At the point $q$, the exponents are $0$ and $2$, so the angle is $2$,
and this point is a removable singularity of $f$ by condition b) after
(\ref{1a}),
so the sides meeting at $q$ are mapped to the same circle.

For the converse statement, suppose that a circular pentagon with $D=H$
is given, with the angles
$(\kappa_1,\kappa_2,
\kappa_3,\kappa_4)$ at $(0,1,x,\infty)$ and $2$ at $q$, such that
the sides meeting at $q$ are mapped by $f$ into the same circle.
Then $f$ extends by reflections
to the universal cover of $\C\backslash\{0,x,1,\infty\}$, and the monodromy
of the extended map is a subgroup of $PSL(2,\C)$, the group of linear-fractional
transformations. This means that the Schwarzian derivative
\begin{equation}\label{schwarz}
R:=\frac{f'''}{f'}-\frac{3}{2}\left(\frac{f''}{f'}\right)^2
\end{equation}
is single valued, and the local behavior at $t_j,\,\infty$ and $q$ implies
that $R$ has poles of order two with
$$R(z)=\frac{1-\kappa_j^2}{2(z-t_j)^2}+\ldots,\quad \mbox{as}\quad z\to t_j,$$
$$R(z)=-\frac{3}{2(z-q)^2}+\ldots,\quad\mbox{as}\quad  z\to q,$$
and similarly at infinity, so $R$ is a rational function. As the intervals
of the real line between the singularities are mapped to arcs of circles,
$R$ is real on the real line. As $t_j,q$ and $\kappa_j$ are real, we conclude
that the residues of $R$ are also real.
Then the general solution of the Schwarz differential equation (\ref{schwarz})
is a ratio of two linearly independent solutions of (\ref{1a}),
see, for example \cite{Golubev}, \cite{HC}.
The condition that the images
of the sides meeting at $q$ belong to the same
circle ensures that $f$ has trivial monodromy at $q$, so $q$ is
an apparent singularity with exponents $0$ and $2$.
This completes the proof of Proposition~1.
\vspace{.1in}

If $C_k$ is the circle containing $f([t_{k-1},t_k])$ then
$(C_1,C_2,C_3,C_4)$ is a chain of
{\em four} circles for which (\ref{C1}) holds and
\begin{equation}
\label{C2}
C_j\neq C_{j+1},\quad j\in\Z_4
\end{equation}
in view of the condition (\ref{as}). If we denote
by $\sigma_j$ the reflection in $C_j$, then
the projective monodromy generators are
\begin{equation}\label{gen}
T_j=\sigma_{j}\sigma_{j+1},\quad j\in\Z_4.
\end{equation}
This assumes that the fundamental group generators are chosen as
in Fig.~\ref{loops}. To prove (\ref{gen}) we notice that each loop
$\gamma_j,\; 1\leq j\leq 3,$ first crosses the interval $(t_{j-1},t_j)$
from $H$ to the lower half-plane, and crosses $(t_j,t_{j+1})$ back to $H$.
The first crossing corresponds to the reflection $\sigma_j$ and the second
to the reflection in the circle $\sigma_jC_{j+1}$. This second reflection
is $\sigma_j\sigma_{j+1}\sigma_j$, so the whole continuation
around $t_j$ is performed with the reflection
$$\sigma_j\sigma_{j+1}\sigma_j\sigma_j=\sigma_j\sigma_{j+1},$$
as stated.

If a projective monodromy representation satisfies (\ref{gen})
with some reflections $\sigma_j$, we say that
this representation is {\em generated by reflections}.
In Appendix I we will find the necessary and sufficient conditions
for $T_j$  to be generated by reflections, and will show how to
find the $\sigma_j$ when these conditions hold. We will see that
the reflections $\sigma_j$ are uniquely defined by the monodromy generators,
except in the case when all these generators commute.
\vspace{.1in}

\noindent
{\bf 5. Special pentagons}
\vspace{.1in}

The previous section motivates
consideration of
pentagons with one angle equal to
$2$, and the sides forming this angle mapped into the same
circle by the developing map, while each pair of sides meeting at
one of the other corners is mapped by the developing map to distinct circles.

We call them {\em special pentagons} and use the following notation
$$Q=(\overline{D},t_1,t_2,t_3,t_4;q,f),$$
where $t_j$ are naturally ordered corners with angles $\alpha_j$,
while $q$ is the corner with angle $2$ which can lie on any arc
between the $t_j$, and the sides meeting at $q$ are mapped by $f$
to the same circle, while the circles containing $f([t_{j-1},t_j])$
and $f([t_j,t_{j+1}])$ are distinct for all $j\in\Z_4$.

This notation is slightly inconsistent with our general
notation for a circular pentagon, because only $t_j$ are
listed in their natural order, while $q$ can be on any
interval between them. To stress this, we separate $q$ from the $t_j$
by a semi-colon.

We recall that a {\em conformal quadrilateral}$\,$\footnote{Not to be confused
with circular quadrilateral!} is a simply connected
Riemann surface, which is conformally equivalent to a disk, with
4 marked prime ends\footnote{For a Jordan region in the plane
prime ends are just boundary points. For general simply connected regions
we refer to \cite{Ahl} and Appendix III.} Conformal equivalence of conformal quadrilaterals
means the existence of a conformal map between them
which maps the marked points
to the marked points.

Each conformal quadrilateral is conformally
equivalent to a rectangle whose marked boundary points are the corners.

We consider special pentagons $Q=(\overline{D},t_1,\ldots,t_4;q,f)$
as conformal quadrilaterals
$(\overline{D},t_1,\ldots,t_4)$,
(forgetting the corner $q$) and define the {\em modulus}
$$\mod Q\in (0,\infty)$$ as the extremal distance in $D$ between the segments
$[t_1,t_2]$
and $[t_3,t_4]$ of $\partial D$.
For the definition and general properties of the extremal distance
we refer to \cite{Ahl}.

To avoid confusion with the sides of a pentagon as defined before,
we use the word {\em segments} to denote
$[t_j,t_{j+1}]\subset\partial D$.
Thus one of the segments consists of two sides of the pentagon and
contains $q$, while each of the
other three segments is the closure of one side of the pentagon.


Every conformal quadrilateral is equivalent to $(\overline{H},0,1,x,\infty)$
for some $x\in(1,\infty)$. With our convention that
$(t_1,t_2,t_3,t_4)=(0,1,x,\infty)$,
the modulus
is a strictly
increasing function of $x$, mapping $(1,+\infty)$ onto $(0,+\infty)$.
An explicit expression of this function can be found in \cite{Ahl} but
we do not need this formula.

To state the properties of the extremal distance that we need,
we use the intrinsic distance on $\overline{D}$
defined in section 3.
\vspace{.1in}

\noindent
{\bf Lemma 1.} (\cite[Lemma 13.1]{EGT1} and Lemma A4 in Appendix III.)
{\em Consider a sequence of special pentagons
$Q_n$ whose areas are bounded from above. If the
intrinsic distance between $[t_1,t_2]$ and $[t_3,t_4]$
tends to zero, while the
intrinsic distance between $[t_2,t_3]$ and $[t_4,t_1]$
stays
away from zero, then $\mod Q_n\to 0.$ If the intrinsic distance
between $[t_2,t_3]$ and $[t_4,t_1]$
tends to $0$  while the intrinsic distance between
$[t_1,t_2]$ and $[t_3,t_4]$
stays away from zero then $\mod Q_n\to \infty$.}
\vspace{.1in}

\noindent
{\bf 6. Evolution of special pentagons. Local families}
\nopagebreak
\vspace{.1in}

\nopagebreak
We recall that $f(q)$ is called the {\em tip} of the slit. The slit
can be lengthened or shortened with $f(q)$ moving
on a circle.
Lengthening or shortening the slit along the circle while keeping
all circles of the chain unchanged,
we obtain a one-parametric family of special pentagons, parametrized by
some interval. We choose the length of the slit as parameter.
\vspace{.1in}

\noindent
{\bf Lemma 2.} {\em As a function of the length
of the slit, $\mod Q$ is monotone.
It is strictly increasing if $q\in (t_2,t_3)\cup(t_4,t_1)$
and strictly decreasing if $q\in (t_1,t_2)\cup(t_3,t_4)$.}
\vspace{.1in}

This follows from the standard properties of the modulus, \cite[4.3]{Ahl}
and Lemma A3 in Appendix III.

As the slit shortens, it eventually vanishes, and
we obtain a polygon with
at most $4$ sides.
As the slit becomes longer, it eventually hits the boundary and becomes
a cross-cut which splits $\overline{D}$ into two polygons.

Such a family, obtained from a special pentagon by shortening the slit
until it vanishes and lengthening the slit until it hits the boundary,
will be called a {\em local family} of special pentagons.
It is parametrized by an open interval (for example, the length of the slit),
and corresponds to an open interval on the ray $(1,+\infty)$
in view of Lemma 2.

In the remainder of this section we will study in detail
what happens at the ends of a local family.
In the next section we will see how local families are combined
into a {\em global family}
of special pentagons, parametrized by $x\in(1,+\infty)$,
so that the special pentagons of the global family depend continuously
and even real-analytically on $x$.

Consider a local family of special pentagons $Q_x$ parametrized by
$x\in I$ where $I$ is an interval
in $(1,\infty)$.

We say that the modulus {\em degenerates} if $\mod Q_x\to 0$
or $\mod Q_x\to\infty$ as $x$ tends to an endpoint of $I$.
This means that this endpoint must be $1$ or $\infty$.

First we state the conditions of degeneracy.

Suppose that $q\in(t_{k-1},t_k)$. Suppose that the slit
shortens and vanishes, then
$q$ must collide with a corner $t_{k-1}$ or $t_k$.
If the intrinsic length of $[t_{k-1},q]$ is strictly smaller than the intrinsic
length of $[q,t_k]$ then $q$ will collide with $t_{k-1}$.
If the intrinsic length of $[t_{k-1},q]$ is strictly greater than the intrinsic
length of $[q,t_k]$ then $q$ will collide with $t_k$.
In both cases we obtain a non-degenerated quadrilateral in the limit,
so $x$ tends to some $x_0\in(1,\infty)$, and $x_0$ is a special point
with the value $q(x_0)=t$, where $t\in\{ t_{k-1},t_k\}$ is the corner
with which $q$ collided. If the intrinsic lengths of
$[t_{k-1},q]$ and $[q,t_k]$ are equal, then as the slit shortens and vanishes,
$t_{k-1}$ and $t_k$ collide, and the limit polygon is a triangle.

The degeneracy condition is thus the following:
\vspace{.1in}

{\em D1. When the slit shortens and vanishes, $\mod Q_x$ degenerates
if and only if two corners collide in the limit}.
\vspace{.1in}

In other words, our special pentagon must be a slit triangle with
the slit originating at some vertex $A$. Notice that the angle of
the triangle at $A$ may be an integer, and the images of sides of the
triangle which are adjacent at $A$ may belong to the same circle.

Now suppose that the slit lengthens. Then eventually it will hit the
boundary from inside at some point $s\in\partial{H}$.
\begin{eqnarray}\label{hits}
&&\text{\em This means that}\; f(s)\;
\text{\em belongs to the circle of the slit}\\
&&\text{\em
and the intrinsic distance between}\; q\; \text{\em and}\; s\; \text{\em tends to zero.}\nonumber
\end{eqnarray}

Suppose that $q\in(t_{k-1},t_k)$. If $s\in[t_{k+1},t_{k+2}]$,
then the modulus degenerates, otherwise it does not.
So we have the second degeneration condition:
\vspace{.1in}

{\em D2. When the slit lengthens and splits the pentagon,
$\mod Q_x$ degenerates if and only if the slit hits
the segment which is opposite to the segment to which $q$ belongs.}
\vspace{.1in}

In other words, in the limit, the slit splits the boundary into two arcs, and the
modulus degenerates if and only if the closures of these two arcs
contain at least two corners each.

Now we consider non-degenerate cases, that is the cases when there exists
a limit quadrilateral when the slit vanishes or when it splits the pentagon.
\vspace{.1in}

{\em Case 1. The slit vanishes, $q$ collides with exactly one corner,
and the angle of the special pentagon at this corner is positive.
}
\vspace{.1in}

{\em Case 2. The slit lengthens and hits the boundary at an interior point
of the side, splitting the special pentagon into a quadrilateral
and a digon with positive angle.}
\vspace{.1in}

It is also possible that as the slit lengthens, it hits the boundary
at a corner. If the modulus does not degenerate, this must be a corner
neighboring $q$, for example $t_k$. Then the special pentagon splits
into a non-degenerate quadrilateral and the remaining part which can
have only one corner at $t_k$. Therefore
the detached part must be a disk.
This we call
\vspace{.1in}

{\em Case 3. The slit lengthens and hits the boundary at a corner. A disk
splits away from the special pentagon, leaving a non-degenerate quadrilateral.}
\vspace{.05in}

The remaining cases happen when the slit vanishes as in Case 1,
and the corner with which the slit collides has zero angle, or
when the slit lengthens, splits the special pentagon,
and one part of the split pentagon is a digon with zero angle.
These cases will be considered in the next section.
\vspace{.1in}

\noindent
{\bf 7. Transformations connecting local families into a global family.}
\vspace{.1in}

In this section we explain what happens when $x$ passes a special point $x_0$,
and $q$ passes one of the $t_j$, so that $Q_{x_0}$ is a non-degenerate
quadrilateral.

We describe four types of transformations that may occur. The first three
correspond to cases 1-3 of the previous section, and the 4-th transformation
to the two remaining cases with a zero angle.
\vspace{.1in}

\begin{figure}
\centering
\includegraphics[width=5in]{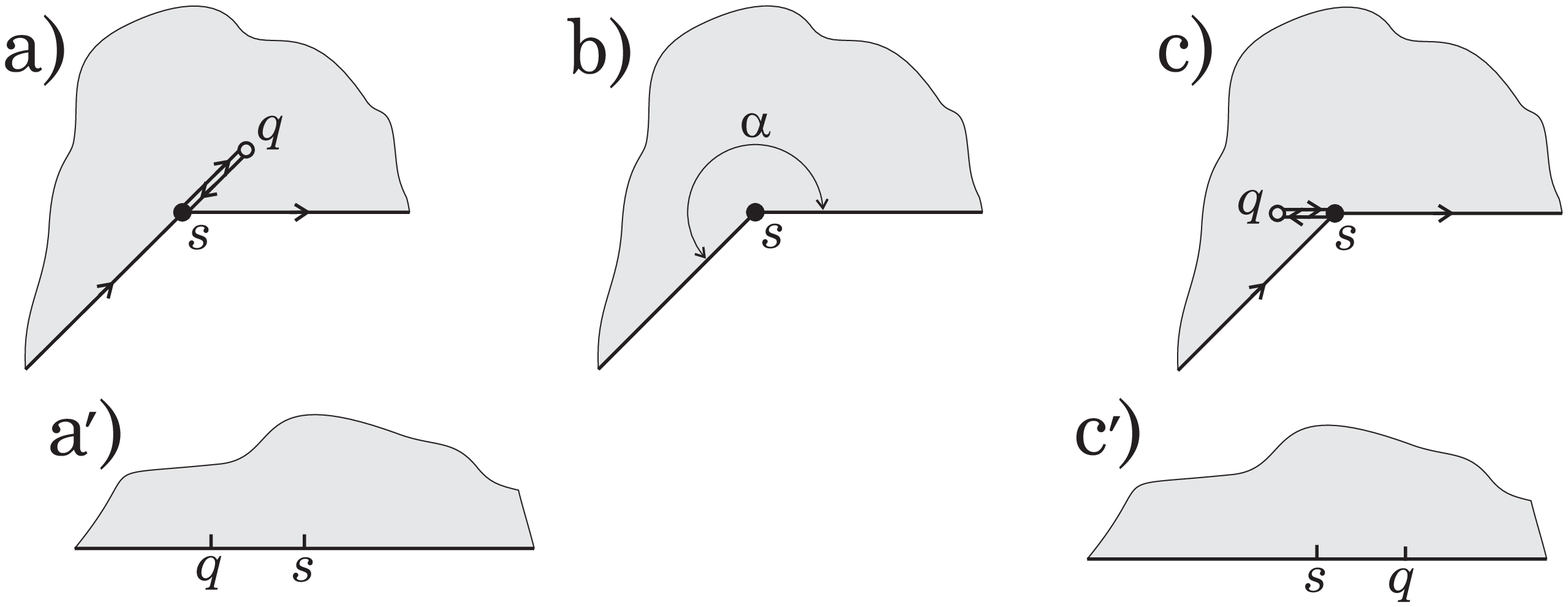}
\caption{Transformation 1.}\label{case1}
\end{figure}

\noindent
{\em Transformation 1 (Fig.~\ref{case1}).  Vanishing slit, Case 1.}
\vspace{.1in}

Suppose that $x$ passes a special point
$x_0$. Before this $q\in(t_{k-1},t_k)$, and when $x=x_0$, $q(x_0)=s$,
where $s$ is one of the points $t_{k-1},\; t_k$.
After $x$ passes $x_0$, $q$ and $s$ interchange. As the images of the sides
lie on the same fixed circles, we have the situation shown in Fig.~\ref{case1}.
\vspace{.1in}

{\em The slit whose image was an arc of a circle vanishes, and then
a new slit starts growing
with the image on an arc of the circle that is adjacent
to the previous circle at the
corner where the old slit vanished.}
\vspace{.1in}

The angle $\alpha$ of the limit quadrilateral
at the corner where the slit vanishes
satisfies
\begin{equation}\label{alphag1}
\alpha>1.
\end{equation}

\begin{figure}
\centering
\includegraphics[width=5in]{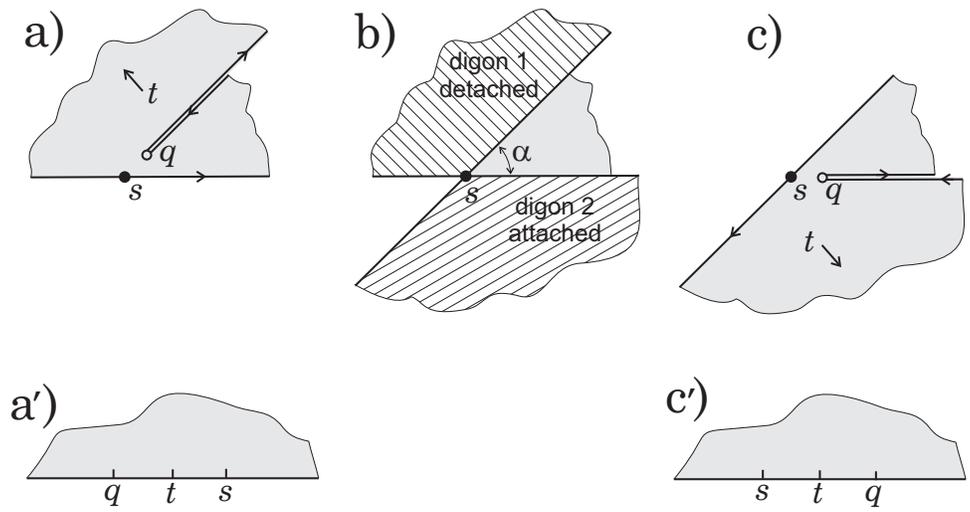}
\caption{Transformation 2.}\label{case2}
\end{figure}

\noindent
{\em Transformation 2 (Fig.~\ref{case2}). The slit hits
an interior point $s$ of a segment (see (\ref{hits})). The slit is not tangent
to this segment (Case 2).}
\vspace{.1in}

If the modulus does not degenerate, $s$ must be an interior point
of the segment
adjacent to that segment which contains $q$.
So there is exactly one corner $t$ in the interior of one boundary
arc $Z$ between $q$ and $s$, and three corners
on the complementary arc.
When $q$ hits $s$, our pentagon splits into two polygons: a quadrilateral
with the corners $s$ and  $t_j\neq t$, and a digon with the corners $s$ and
and $t$. Notice that the angles at the two corners of a digon are
always equal. It is clear that in the described situation this digon angle is
less than $1$: it is the inclination of the slit to the side that it hits.
Thus if $\alpha$ is the angle of the limit quadrilateral at $s$
then
\begin{equation}\label{alphal1}
\alpha<1,
\end{equation}
while the digon has two angles $1-\alpha$, at $s$ and at $t$.
\vspace{.1in}

{\em When the slit hits the boundary at an interior point of the side,
a digon is detached, and a vertical digon\footnote{
A pair of digons with equal angles formed by two circles are
called {\em vertical}.}
 is attached on the
side which was hit. One side of the old slit becomes a side of
the new pentagon.}

\vspace{.1in}

\begin{figure}
\centering
\includegraphics[width=5in]{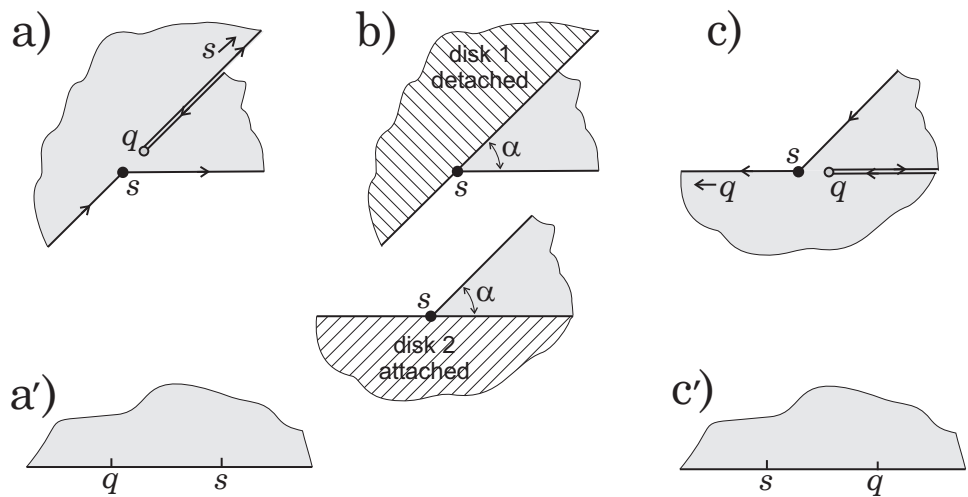}
\caption{Transformation 3.}\label{case3}
\end{figure}

\noindent
{\em Transformation 3 (Fig.~\ref{case3}).
The  slit hits  a corner (see (\ref{hits})) as described in Case 3.}
\vspace{.1in}

In this case $s$ is a corner, and there is no other corner on $Z$.
When $q$ hits $s$, the special pentagon splits into
a quadrilateral with corners at $t_j,\;1\leq j\leq 4$ and the other
part which must be a disk.
So if the limit quadrilateral has angle $\alpha$ at $s$ then
the special pentagon before the limit has angle $\alpha+1$
at $s$.

So far we ignored the non-generic cases which
may occur when some circles of the chain
are tangent: when the slit vanishes at a corner with zero angle,
and when the slit hits from inside a side which is tangent to it.
In these cases one more transformation occurs.
\vspace{.1in}

\begin{figure}
\centering
\includegraphics[width=5in]{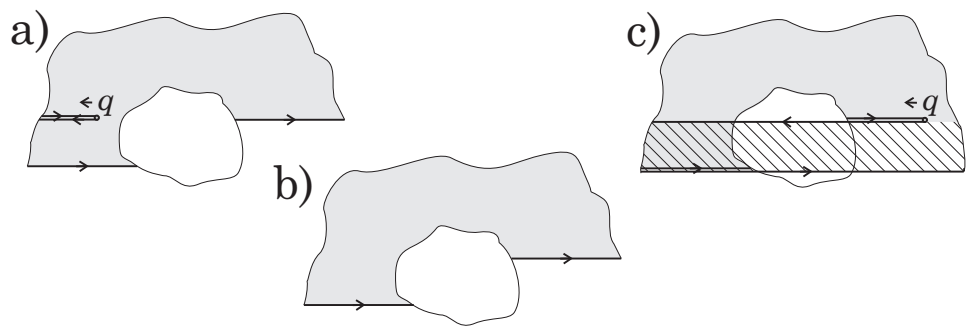}
\caption{Transformation 4.}\label{case4}
\end{figure}

\noindent
{\em Transformation 4 (Fig.~\ref{case4}). The slit vanishes at a corner with zero angle
and a digon with zero angle is attached.}
\vspace{.1in}

The slit shortens and vanishes at the corner with zero angle
(which is shown at $\infty$
in Fig.~\ref{case4}a,
and the resulting quadrilateral in Fig.~\ref{case4}b has angle $1$.
After that, we attach to this quadrilateral a digon with zero angle
(shown as a strip in Fig.~\ref{case4}c,
and the slit shortens when $x$ continues to change in the same direction.

When we run $x$ backwards, we first encounter Fig.~\ref{case4}c with
the lengthening slit which hits the boundary of the pentagon from inside
under zero angle. Similarly to transformation 2, a digon detaches
(the strip in Fig.~\ref{case4}c is a digon with zero angle),
and a new slit starts growing as in Fig.~\ref{case4}a.

Notice that unlike in all other transformations 1-3, the direction of
the slit evolution (whether it lengthens or shortens) does not change
for this transformation. Two other distinctions of this transformation
from transformations 1--3 are that the old and new slit are on the same
circle, and that $q$ is on the same segment before and after
the transformation.

This is consistent with the fact that the
special points of the function $q(x)$ are simple, unless the
angle of the special pentagon corresponding to a special point is zero,
in which case this special point is a simple critical point for
$q(x)$ \cite[Ch. 9, \S 46]{Gromak}.
The process we described shows that every local family
of special pentagons can be extended
to a global family of special pentagons,
with the special pentagons becoming
quadrilaterals at isolated
points. At these points
one angle $\alpha$ of the pentagon becomes angle $\beta$ of
the quadrilateral, and these angles are related as follows:
$\beta=\alpha+1$
for transformation 1,
$\beta=\alpha-1$ for transformation 3, and
$\beta=1-\alpha$ for transformation 2.

This continuation can be either performed indefinitely
in one or both directions,
or the modulus can degenerate at one or both ends.

In sections 8--10 we will analyse global families.
\vspace{.1in}

\noindent
{\em Remark.} Transformations 1--4 suggest the following:
\vspace{.1in}

{\em When equation (\ref{1a}) undergoes an isomonodromic deformation and
$q$ collides with some $t_k$, then the resulting Heun equation
has exponent difference $\pm(\kappa_k+1)$ at $t_k$ as in Case 1,
or $\pm(\kappa_k-1)$ as in Cases 2,3}.
\vspace{.1in}

This is true in general, without
our restriction that the $\kappa_j$ and $x$ are real. To obtain
this result one can use asymptotics of $p(x)$ and $q(x)$ as $x$
tends to a special point written in \cite[p. 534-535]{Okamoto2},
and obtain the limit equation with four singularities
directly from (\ref{1a}).
\vspace{.1in}

\noindent
{\bf 8. Explicit description of global families
}
\nopagebreak
\vspace{.1in}

\nopagebreak
The previous section shows how local families are combined
into a global family. A global family $Q(x),\; 1<x<+\infty$ consists
of local families of special pentagons
parametrized by intervals $x_j<x<x_{j+1}$.
At the points $x_j\in(1,\infty)$
the special pentagon becomes a non-degenerate quadrilateral.
A global family may consist of a single local family; such families
will be discussed
in section 10.
The sequence $x_j$ can be finite or infinite
in one direction or in both directions. The smallest and the largest
terms of this sequence, when they exist, are $1$ and $\infty$. All other
terms correspond to the special points of the solution
of PVI which is described
by our global family.
To describe global families more precisely, we recall
the construction of combinatorial objects related to circular polygons.
\vspace{.1in}

\noindent
{\bf 9. Representation of polygons by nets.}
\vspace{.1in}

Circular polygons are conveniently represented by nets
\cite{EGA,EGT1,EGT2}.
Consider a polygon given by (\ref{ngon}). Its {\em net} is the cell
decomposition of $\overline{D}$ by all $f$-preimages of the
circles $C_k$, where $C_k$ is the circle that contains
$f([t_{k-1},t_k])$. The corners are required to be vertices
of the net. The $1$-cells of the net are labeled by their images.
Two nets are considered equal if there is an orientation preserving
homeomorphism which maps one onto another preserving the labels of
the corners. Let $e_1$ be the $1$-cell on the boundary of the
net, oriented according to the orientation of the boundary,
and beginning at $t_1$.
Two polygons with developing maps
$f_1,f_2$ are equal if their circles $C_j$ are the same,
their nets are equal, and the images $e_1$, as oriented $1$-cells,
are equal.

In the illustrations we label
the circles $C_j$ and corresponding $1$-cells of the nets with colors
(or with different styles of lines in the black and white version).

\begin{figure}
\centering
\includegraphics[width=3.5in]{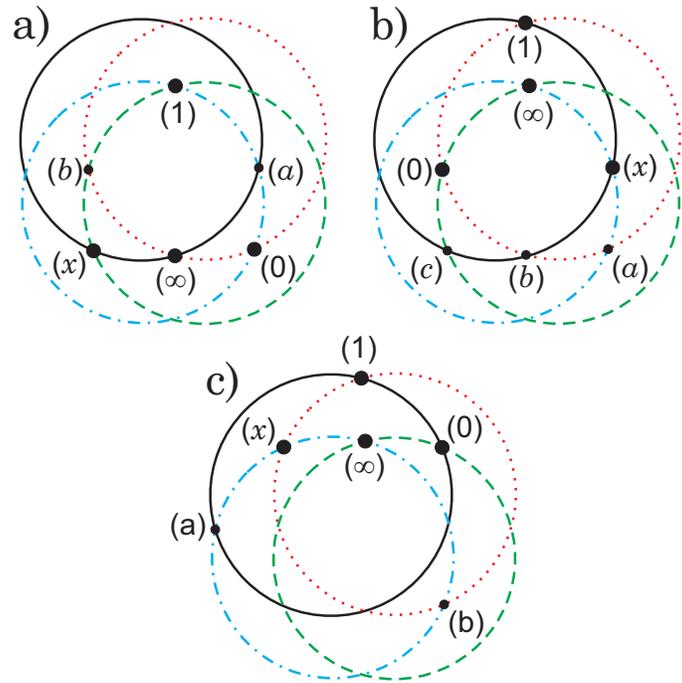}
\caption{Four-circle chains for a) Example 1, b) Example 2, c) Example 7.}\label{4circles-pent}
\end{figure}

It is difficult to characterize all possible nets of circular quadrilaterals
or special pentagons. The topological
classification of generic $4$-circle chains
is given in Appendix II. For each type of chain, one has a set of nets
compatible with this chain. For the chain topologically equivalent to
a quadruple of generic great circles as in Fig.~\ref{4circles-pent},
one can give the following characterization
of the nets. Notice that the cell decomposition of the sphere
in Fig.~\ref{4circles-pent} has the following property which must be inherited
by the net:
\vspace{.1in}

a) When two $2$-cells share a boundary $1$-cell, one of these two
cells is a quadrilateral and another is a triangle.
\vspace{.1in}

Moreover, the net has an evident additional property:
\vspace{.1in}

b) All interior vertices have degree $4$, and all vertices on the
sides have degree $3$.
\vspace{.1in}

Our standing assumption (\ref{as}) implies that
\vspace{.1in}

c) The degrees of the corners (as vertices of the net) are even.
\vspace{.1in}

\begin{figure}[H]
\centering
\includegraphics[width=5.1in]{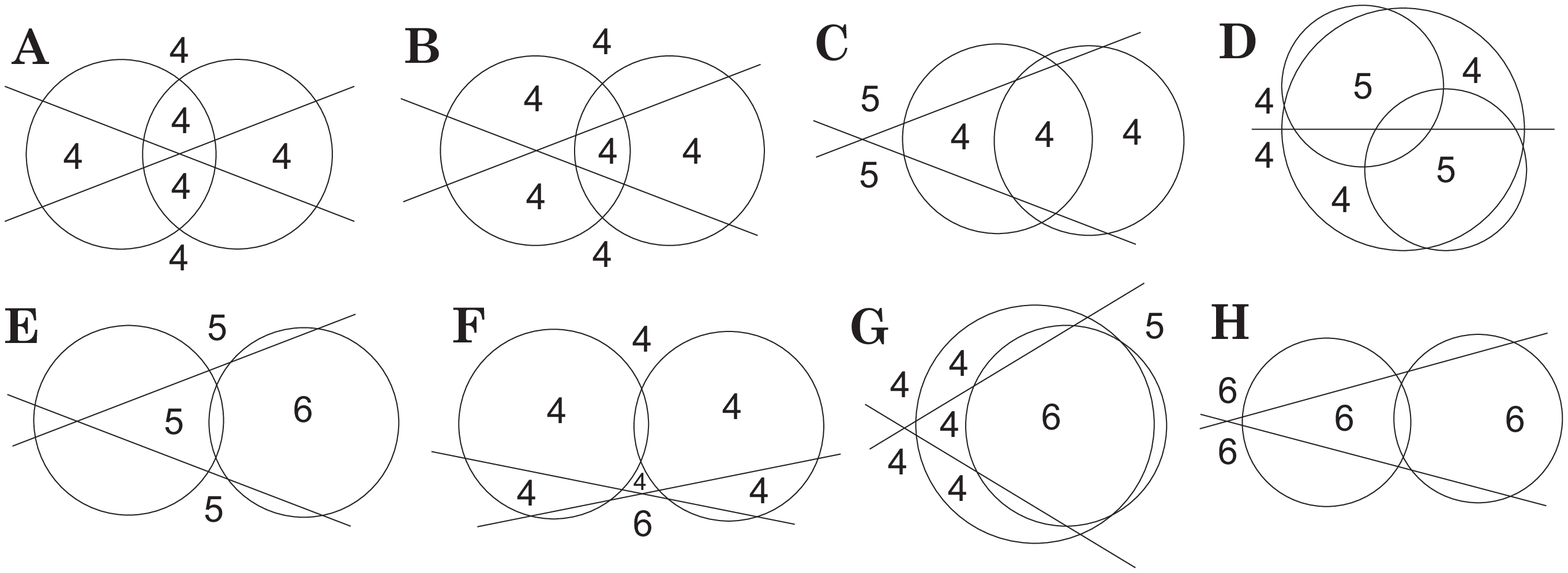}
\caption{Generic chains A-H. All pairs of circles intersect.}\label{chains1}
\end{figure}

These three properties completely characterize the nets of circular
quadrilaterals
over the $4$-circle chain
shown in Fig.~\ref{4circles-pent}. See also Fig.~\ref{chains1}A,
where the same $4$-circle
chain is shown.

\begin{figure}
\centering
\includegraphics[width=5in]{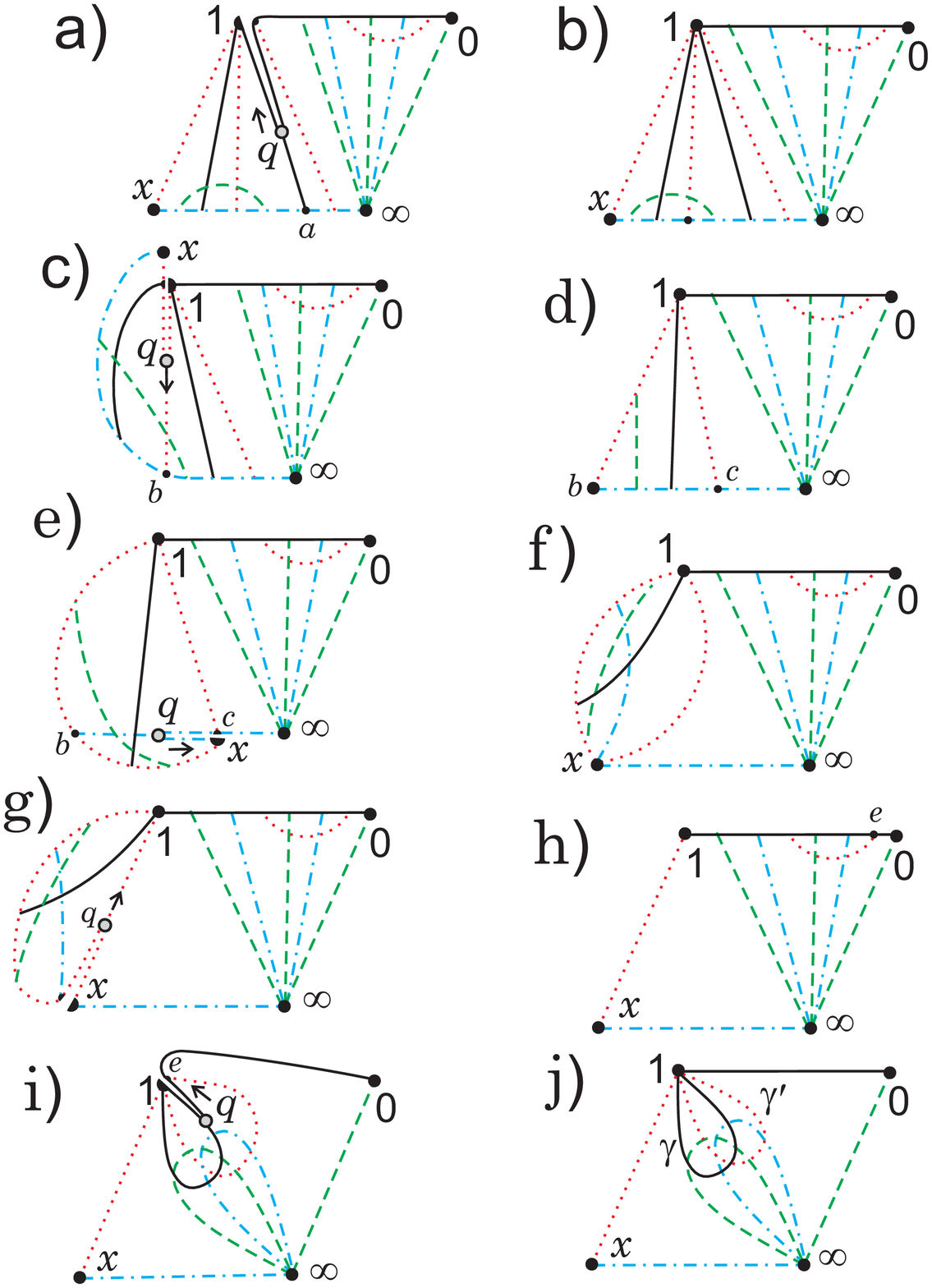}
\caption{Nets for a half of the global family in Example 7.}\label{img15}
\end{figure}

Thus for example, all cell decompositions
in the right column of Fig.~\ref{img15} are nets of quadrilaterals.
\vspace{.1in}

\noindent
{\bf 10. Real solutions of PVI without special points}
\vspace{.1in}

Our paper \cite{EGH} describes
all complex solutions of PVI without special points in the complex plane.

In this section we will describe
all real solutions of PVI with real parameters, which
have no real special points.

For simplicity we limit ourselves to the generic case: all parameters
$\kappa_j$ are not integers, and the circles of the chain are not tangent
to each other. The last condition holds for example when the projective
monodromy contains no parabolic transformations.

Solutions without special points
correspond to local families for which the modulus
degenerates on both ends. So degeneracy conditions D1 and D2
of section 6 must be satisfied (one condition on one end
and another on another end).

\begin{figure}
\centering
\includegraphics[width=5in]{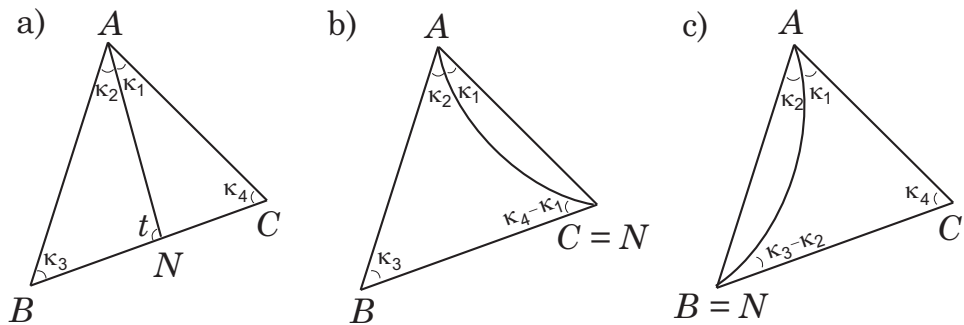}
\caption{Slit triangles.}\label{triang}
\end{figure}

Thus we have one of the three configurations shown in
Fig.~\ref{triang}.

In a family without special points we must have
$$q(x)\in (t_j,t_{j+1})\quad\mbox{for all}\quad x\in(1,\infty)$$
and some $j\in\Z_4$. Suppose without loss of generality that
$j=1$, so that
\begin{equation}
\label{between}
q(x)\in(0,1),\quad x>1.
\end{equation}
The condition that $f(0)=f(1)$ when the slit vanishes
can be stated in terms of projective monodromy representation:
\begin{equation}\label{00}
\begin{split}
\mbox{Transformations}\; T_0\;\mbox{and}\;T_1\;\mbox{have a common fixed point}\\
\mbox{with multipliers}\; e^{2\pi i\kappa_1}\;\mbox{and}\;e^{2\pi i\kappa_2}.
\end{split}
\end{equation}

We will use the results of Klein \cite{Klein} and Van Vleck \cite{VV}
on circular triangles. First of all we have
\vspace{.1in}

\noindent
{\bf Lemma 3.} {\em For any positive numbers $\lambda_j,\; 1\leq j\leq 3,$
there exists a unique equivalence class of circular triangles
with angles $(\lambda_1,\lambda_2,\lambda_3).$}
\vspace{.1in}

{\em Sketch of the proof.} The developing map of a triangle with
angles $\lambda_j$ satisfies the Schwarz differential equation:
$$\frac{f'''}{f'}-\frac{3}{2}\left(\frac{f''}{f'}\right)^2=
\frac{1-\lambda_1^2}{2z^2}+\frac{1-\lambda_2^2}{2(z-1)^2}+
\frac{\lambda_1^2+\lambda_2^2-\lambda_3^2-1}{2z(z-1)},$$
see, for example, \cite[p. 452]{HC} or \cite[Ch. VI, \S 3]{Golubev}.
Parameters $\lambda_j$ can be arbitrary non-negative numbers,
and for fixed parameters all solutions give equivalent triangles.
This proves the
lemma.
\vspace{.1in}

\noindent
{\bf Lemma 4.} {\em In a triangle with
angles $(\lambda_1,\lambda_2,\lambda_3)$, the image of the closure of
the side opposite to $\lambda_1$ under the developing map makes
$$E\left(\frac{\lambda_1-\lambda_2-\lambda_2+1}{2}\right)$$
full turns\footnote{
When $f:[a,b]\to C$ is
an immersion of a closed interval $[a,b]$ into a circle $C$
\def\card{\mathrm{card}\, }
then the ``number of full turns'' is defined as $\card f^{-1}(a)-1$.}
around the circle containing this image.
Here $E(x)$ is the integer part of $x$ for $x>0$ and zero otherwise.}
\vspace{.1in}

This is called the Erg\"anzungsrelation of Klein \cite{Klein}, \cite{VV}.

First we address the easier cases b) and c) in Fig.~\ref{triang}. Consider the
case b). In this case,
our triangle $ABC$ is split into two parts, one of which is a digon with
both angles $\kappa_1$,
and the other part is a triangle with angles $(\kappa_2,\kappa_3,
\kappa_4-\kappa_1)$.
The image of a side of a digon cannot cover a full circle. Therefore
the image of the side $NA$ of the triangle cannot cover the full circle.
The necessary and sufficient condition for this according to Lemma~4
is
\begin{equation}\label{33}
\kappa_3+\kappa_1\leq\kappa_2+\kappa_4+1.
\end{equation}
Similarly, Fig.~\ref{triang}c produces the condition
\begin{equation}\label{44}
\kappa_4+\kappa_2\leq \kappa_1+\kappa_3+1.
\end{equation}
As there are no free parameters in the configurations in Fig.~\ref{triang}b,c,
we conclude that when (\ref{33}) is satisfied, there is a single equivalence
class of
configurations of the form Fig.~\ref{triang}b and when (\ref{44})
is satisfied, there is
a single equivalence class of configurations
of the form Fig.~\ref{triang}c with these angles.
This means that each PVI with such parameters has an isolated
solution of type b) or c), or both.
Notice that the two inequalities (\ref{33}) and (\ref{44})
cover the whole range
of real parameters, so we conclude that isolated solutions of PVI of
one or both types Fig.~\ref{triang}b,c always exist.
There can be one or two of them.
\vspace{.1in}

Now we turn to the case a).
We introduce the auxiliary angle $t\in[0,1]$ as a parameter
(see Fig.~\ref{triang}a).

When $t$ is fixed, there exist two equivalence classes of
triangles with prescribed
angles: $NAB$ with angles $(t,\kappa_2,\kappa_3)$
and $ANC$ with angles $(\kappa_1,1-t,\kappa_4)$.
Let $f_1$ and $f_2$ be their developing
maps.

The question is when we can glue such two triangles along the side $NA$.
As we assume that the circles of the chain are not tangent to
each other,
$t\in(0,1)$, and we can post-compose the $f_j$ with
linear-fractional transformations to achieve
$f_j(N)=0$,
$f_1(BN)$ and $f_2(NC)$ belong to a line $\ell$ through the origin,
and $f_j(NA)$ is contained in the real line for $j=1,2.$
Then it is clear that the necessary
and sufficient condition for the possibility of gluing is that
the image of the closed side $NA$ under both developing maps
intersects the line $\ell$ the same number of times. Indeed, if this is so,
the images of the point $A$ under both developing maps are on the same side
of the line $\ell$  (or both at $0$, or both at $\infty$),
and these images can be made equal by
additional scaling $z\mapsto rz$ with some $r>0$.

As $\ell$ intersects the real line at $0$ and $\infty$ we are interested
in the combined numbers of zeros and poles of $f_j$ on $NA$.

Consider triangle $NAB$ and denote by $(\lambda_1,\lambda_2,\lambda_3)$
the angles at $(N,A,B)$.
To count the number of zeros and poles of $f_1$ on $NA$ times the image of $AN$ intersects the circle
that contains $NB$ we use the results of Klein \cite{Klein} (see also
\cite{VV}) on the number of zeros of hypergeometric function
on an interval $[0,1]$.

We recall that the hypergeometric function $w(z)=F(\alpha,\beta,\gamma,z)$
is the solution of the hypergeometric equation
\begin{equation}\label{hyper}
z(z-1)\frac{d^2w}{dz^2}-(\gamma-(\alpha+\beta+1)z)\frac{dw}{dz}+\alpha\beta w=0,
\end{equation}
which satisfies $F(0)=1$ and is holomorphic at $0$.
A second linearly independent solution of the same equation is
$$F_1(z)=z^{1-\gamma}F(\alpha-\gamma+1,\beta-\gamma+1,2-\gamma,z).$$
Thus $F/F_1$ is a developing map of a triangle whose angles are
the absolute values of the exponent differences of (\ref{hyper}),
$$\lambda_1=|1-\gamma|,\quad \lambda_2=|\gamma-\alpha-\beta|,\quad
\lambda_3=|\alpha-\beta|.$$
We choose
$$\lambda_1=1-\gamma,\quad \lambda_2=\alpha+\beta-\gamma,\quad\lambda_3=\alpha-\beta,$$
which defines $\alpha,\beta,\gamma$ uniquely. Notice that $F_1(0)=0$
because $\gamma\in(0,1)$. Then $f_1=F/F_1$
is the developing map of $NAB$, and the side $NA=(0,1)$.

The number of crossings between $f_1(AN)$ and $\ell$
is
equal to the combined number of zeros and poles of $F$ and $F_1$
on $[0,1]$. We only consider the case $\lambda_1\in(0,1)$ which we need.

According to \cite{VV} the number of zeros of $F$
on $[0,1]$ is:
\vspace{.1in}

\noindent
(i) zero, if $[\lambda_2]>[\lambda_3]$,
\vspace{.1in}

\noindent
(ii) $E\left((\lambda_1+\lambda_3-\lambda_2+1)/2\right)$, if
$[\lambda_3]>[\lambda_2],$
\vspace{.1in}

\noindent
(iii) $0$ or $1$ depending on whether
$E\left((\lambda_1+\lambda_3-\lambda_2+1)/2\right)$
is even or odd, if $[\lambda_2]=[\lambda_3]$.
\vspace{.1in}

The number of zeros of $F_1$ is always
$$E\left(\frac{\lambda_3-\lambda_1-\lambda_2+1}{2}\right)+1.$$

Adding these together we obtain that the number of crossings between
the image of $AN$ and the circle of $BN$ equals to:
\vspace{.1in}

\noindent
$1$, if $[\lambda_2]\geq[\lambda_3]$,
\vspace{.1in}

and to
$$E\left(\frac{\lambda_1+\lambda_3-\lambda_2+1}{2}\right)+
E\left(\frac{\lambda_3-\lambda_1-\lambda_2+1}{2}\right)+1,\quad
\mbox{if}\quad [\lambda_2]<[\lambda_3].$$
\vspace{.1in}

Applying this result to $NAB$ and the similar result to $CAN$, (or rather
to its mirror image), we obtain that
the gluing is possible if and only if one of the following two
conditions holds:
\vspace{.1in}

\noindent
1) $\kappa_2\geq \kappa_3$ and $\kappa_1\geq\kappa_4$,
\vspace{.1in}

\noindent
2) $\kappa_2<\kappa_3$ or $\kappa_1<\kappa_4$, and
\begin{equation}\label{E1}
E\left(\frac{\kappa_3-\kappa_2+t+1}{2}\right)=
E\left(\frac{\kappa_4-\kappa_1-t+2}{2}\right),
\end{equation}
and
\begin{equation}\label{E2}
E\left(\frac{\kappa_3-\kappa_2-t+1}{2}\right)=
E\left(\frac{\kappa_4-\kappa_1+t}{2}
\right),
\end{equation}
\vspace{.1in}

To simplify these conditions, we put $$t=u+1/2,\quad -1/2\leq u\leq 1/2,$$
$$a=\kappa_3-\kappa_2+1/2,\quad b=\kappa_4-\kappa_1+1/2.$$
Then conditions 1)-2) become
\vspace{.1in}

\noindent
$1'$) $a\leq1/2,\; b\leq 1/2$, or
\vspace{.1in}

\noindent
$2'$) $E((a+1+u)/2)=E((b+1-u)/2),$ and $E((a-u)/2)=E((b+u)/2)$.
\vspace{.1in}

\begin{figure}
\centering
\includegraphics[width=3in]{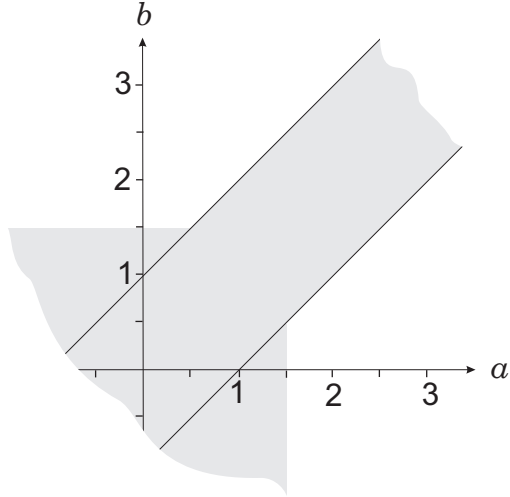}
\caption{Conditions on $a$ and $b$.}\label{non-stair}
\end{figure}

Eliminating $u\in(-1/2,1/2)$ we obtain:
\begin{equation}
\label{condit}
|a-b|<1,\quad\mbox{or}\quad \min\{ a,b\}<3/2,
\end{equation}
the shaded region in Fig.~\ref{non-stair}. For these values of
parameters, PVI has an interval of solutions without special points,
in addition to one or two isolated solutions of types b), c).
The boundaries of the regions corresponding to (\ref{33}) and (\ref{44})
are shown as lines $a=b+1$ and $b=a+1$ in Fig.~\ref{non-stair}.
\vspace{.1in}

\noindent
{\bf Theorem 1.} {\em A real solution $q(x)$ of PVI
with real parameters $\kappa_j$ defined on $(1,\infty)$
and satisfying $q(x)\in (0,1), \; x>1$ always exists.
The monodromy corresponding to this solution satisfies
(\ref{00}). If (\ref{condit}) with $a$ and $b$ as in (\ref{a}), (\ref{b})
holds, then there is an interval of
such solutions.  If
(\ref{33}) or (\ref{44}) hold, then
we have additional one or two isolated solutions.}
\vspace{.1in}

The cases $q(x)\in(t_j,t_{j+1})$ for $j=2,3,4$ are obtained by a
cyclic permutation of the angles in our conditions.
\vspace{.1in}

In \cite{chen}, the following theorem is proved:
{\em Solutions of PVI with parameters $(1/2,1/2,1/2,1/2)$ corresponding
to unitary monodromy do not have special points on a real
interval between the fixed singularities.}

These authors do not assume a priori that their solutions are
real. For the case of real solutions, this result can be obtained as follows.

Consider some special pentagon corresponding to a real solution
of this equation.  The angles are $(1/2,1/2,1/2,1/2,2)$, and
the circles $C_j$ are great circles.
Removing the slit we would obtain a geodesic quadrilateral
with angles $(1/2,1/2,1/2,3/2)$, or a triangle with
angles $(1/2,1/2,1)$, but it is easy to see
that such quadrilateral does not exist (see for example, \cite{EG2}).
So it must be a geodesic triangle with angles
$(1/2,1/2,1)$. Then the special pentagon must have the shape
as in Fig.~\ref{triang}a, so our global family does not
have special points.

In \cite{chen2}, the following fact is proved: {\em Solutions
of PVI with parameters $(1/2,1/2,1/2,3/2)$ corresponding to unitary
monodromy do not have poles on $\R\backslash\{0,1\}$.}
Again, in the case of real solutions, this follows from our results.
When the special pentagon corresponding to such a solution undergoes
a transformation with $q=\infty$ the limit quadrilateral
must have angles $(1/2,1/2,1/2,1/2)$ or $(1/2,1/2,1/2,5/2)$.
Geodesic quadrilaterals with such angles do not exist \cite{EG2}.
\vspace{.1in}

\noindent
{\bf 11. Some other special cases}
\vspace{.1in}

We mention several special cases without going into detail.

1. Suppose that the projective monodromy representation is reducible.
This means that all linear-fractional transformations $T_j$ have a common fixed
point. Without loss of generality, we place it at infinity. Then all
monodromy transformations are affine, and the circles $C_j$
of the chain are straight
lines. Our pentagons are rectilinear and the developing maps can be
expressed by the Schwarz--Christoffel formula. Then $q(x)$ can be expressed
in terms of hypergeometric integrals.

Indeed, the Schwarz--Christoffel formula of a rectilinear special
pentagon gives
$$f(z)=c\int_0^z
\zeta^{\alpha_1-1}(\zeta-1)^{\alpha_2-1}(\zeta-x)^{\alpha_3-1}(\zeta-q)dz=
c(I_1(z)-qI_2(z)),$$
where
$$I_1(z)=\int_0^z
\zeta^{\alpha_1}(\zeta-1)^{\alpha_2-1}(\zeta-x)^{\alpha_3-1}d\zeta$$
and
$$I_2(z)=\int_0^z
\zeta^{\alpha_1-1}(\zeta-1)^{\alpha_2-1}(\zeta-x)^{\alpha_3-1}d\zeta.$$
Normalization $f(1)=a$ and $f(\infty)=b$ will define our special pentagon
completely,
so we obtain with $k=b/a$
$$q(x)=\frac{kI_1(1)-I_1(\infty)}{kI_2(1)-I_2(\infty)},$$
an expression for $q(x)$ in the form of hypergeometric integrals.
\vspace{.1in}

2. Suppose that one of the monodromy transformations is the identity,
and the corresponding angle is $1$.
Then our special pentagon is in fact
a slit triangle (or a slit digon, or a slit disk).
In this case, the developing map itself
can be expressed in terms of hypergeometric functions, see \cite{Schilling},
where the case of slit-triangle quadrilaterals has been studied
in great detail.

These are the two known cases of reduction of PVI when some solutions can be
explicitly found. In a certain sense there are no other cases \cite{W},
except some cases \cite{M2}
when a solution can be expressed as a somewhat non-standard
combination of classical special functions.
\vspace{.1in}

\noindent
{\bf 12. Examples}
\vspace{.1in}



\begin{figure}
\centering
\includegraphics[width=4in]{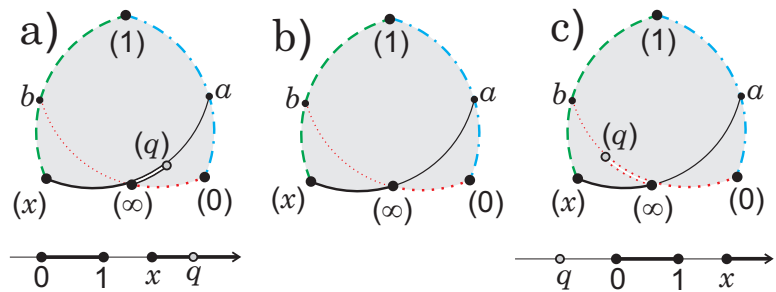}
\caption{Global family in Example 1.}\label{pent1113}
\end{figure}


We begin with the simplest examples when all special pentagons
in a global family are regions in the sphere, so the nets are not
required for their description.
\vspace{.1in}


\noindent{\bf Example 1.}
Consider the conformal map  $f$ from the upper half-plane onto the shaded
region in Fig.~\ref{pent1113}a.
This is the developing map of a special pentagon.
The corners on $\partial H$ are shown below, and their images are shown in
parentheses.

Suppose that the slit lengthens. Then the extremal distance between
the segment $[(0),(1)]$ and the opposite segment $[(x),(\infty)]$
(which contains $(q)$) decreases, so the extremal distance in $H$
between $[0,1]$ and $[x,\infty]$ decreases, thus $x$ decreases.

When $(q)$ tends to $[(0),(1)]$, modulus degenerates and $x\to 1$.

When the slit shortens, $x$ increases, and eventually the slit vanishes
and we obtain Fig.~\ref{pent1113}b. At this moment $q=\infty$.
We have transformation 1, so after that a new slit starts as shown in
Fig.~\ref{pent1113}c and when it hits $[(1),(x)]$, the extremal
distance between $[0,1]$ and $[x,\infty]$
tends to infinity which implies
that $x\to+\infty$.

We conclude that solution $q(x)$ of PVI which corresponds to this
global family has one pole on $(1,\infty)$, and has no zeros,
no fixed points and no $1$-points.

The chain of circles corresponding to this example is shown
in Fig.~\ref{4circles-pent}a. For example, this can be any four
generic great circles, which corresponds to $SU(2)$ monodromy,
determined by formula (\ref{gen}), see Appendix III.

\begin{figure}
\centering
\includegraphics[width=4.5in]{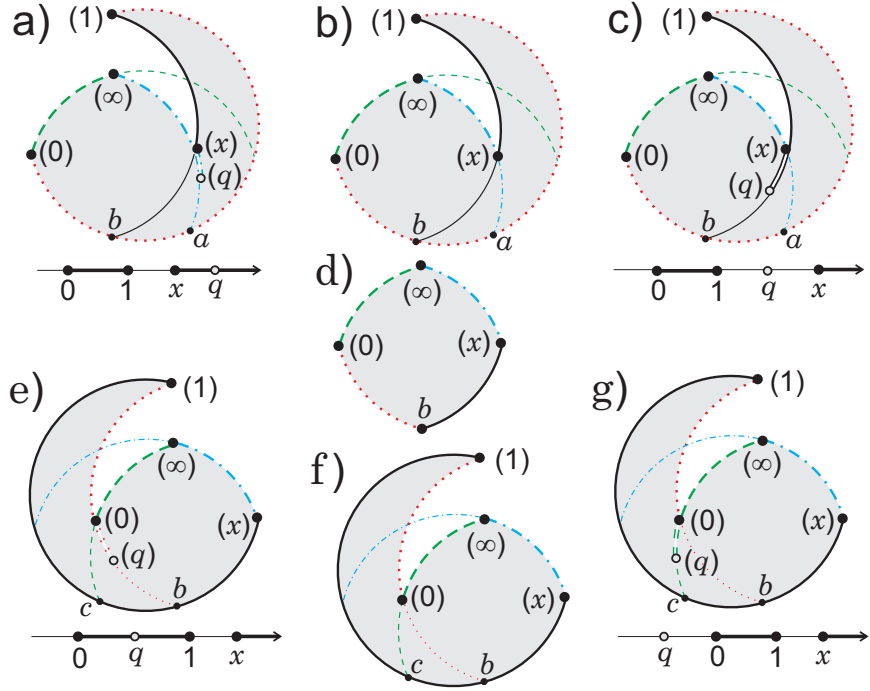}
\caption{Global family in Example 2.}\label{pent1111}
\end{figure}

\begin{figure}
\centering
\includegraphics[width=2.in]{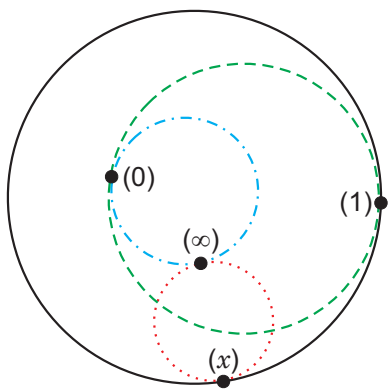}
\caption{Circles chain for Example 3.}\label{zeroangles1}
\end{figure}
\begin{figure}
\centering
\includegraphics[width=4.5in]{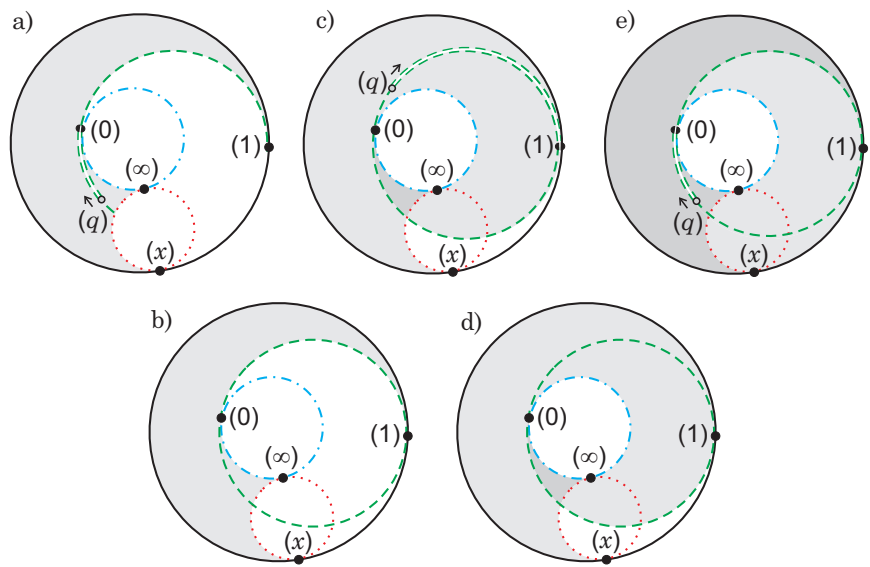}
\caption{Example 3.}\label{zeroangles2}
\end{figure}


\begin{figure}
\centering
\includegraphics[width=2in]{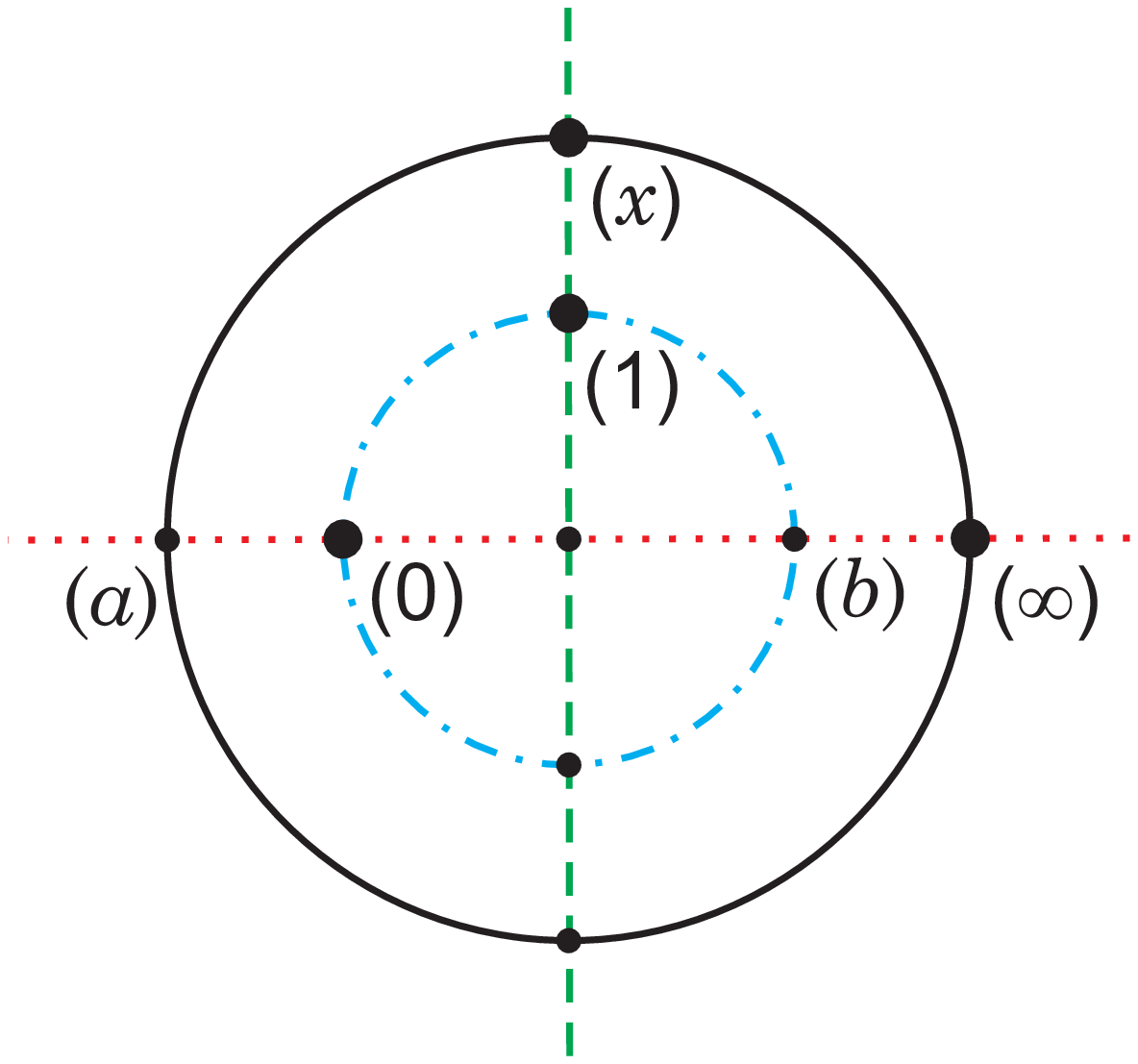}
\caption{Four-circle chain for Example 4.}\label{4circles-pent-b}
\end{figure}

\begin{figure}
\centering
\includegraphics[width=4in]{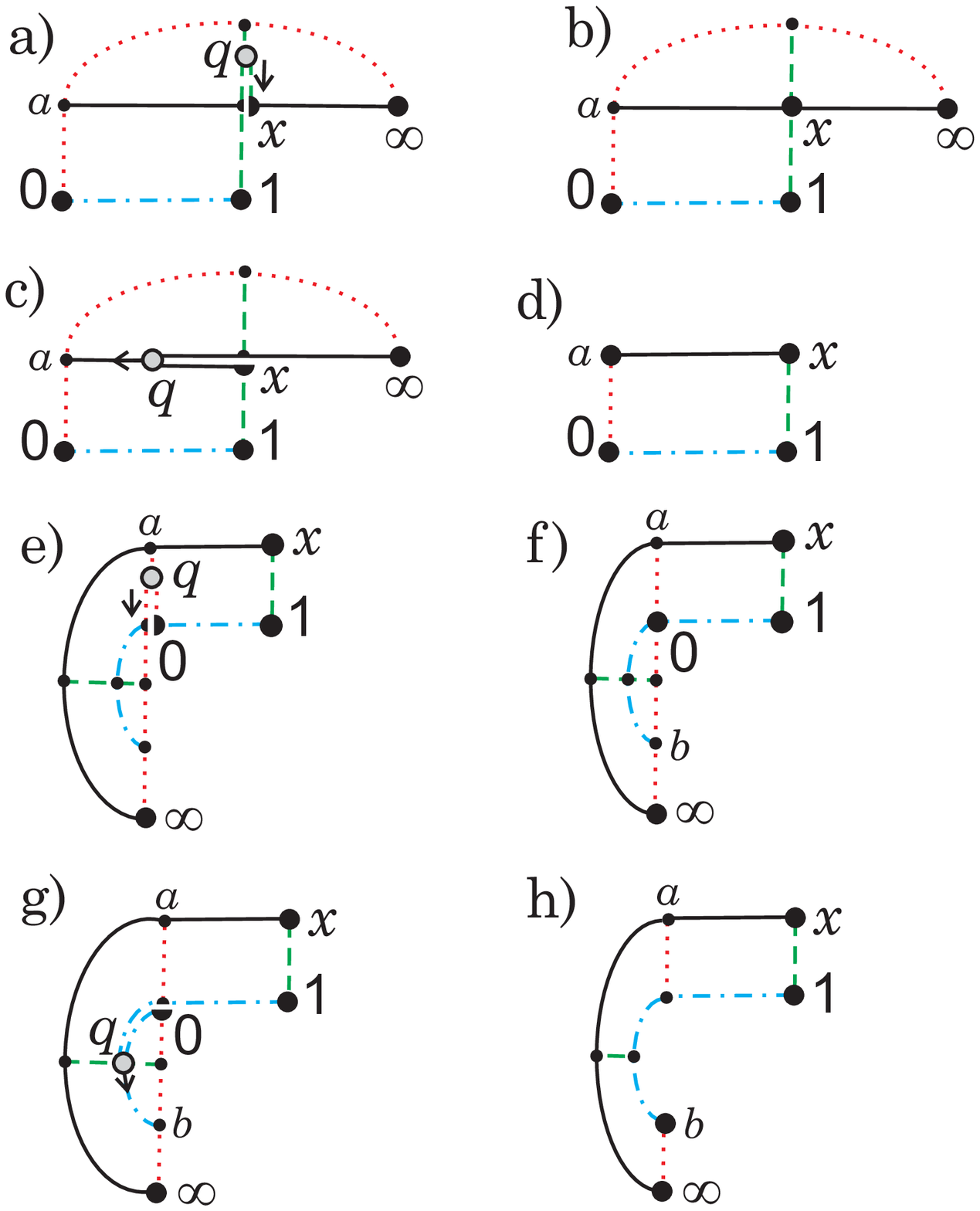}
\caption{Nets for the global family in Example 4.}\label{img12-b}
\end{figure}

\begin{figure}
\centering
\includegraphics[width=3in]{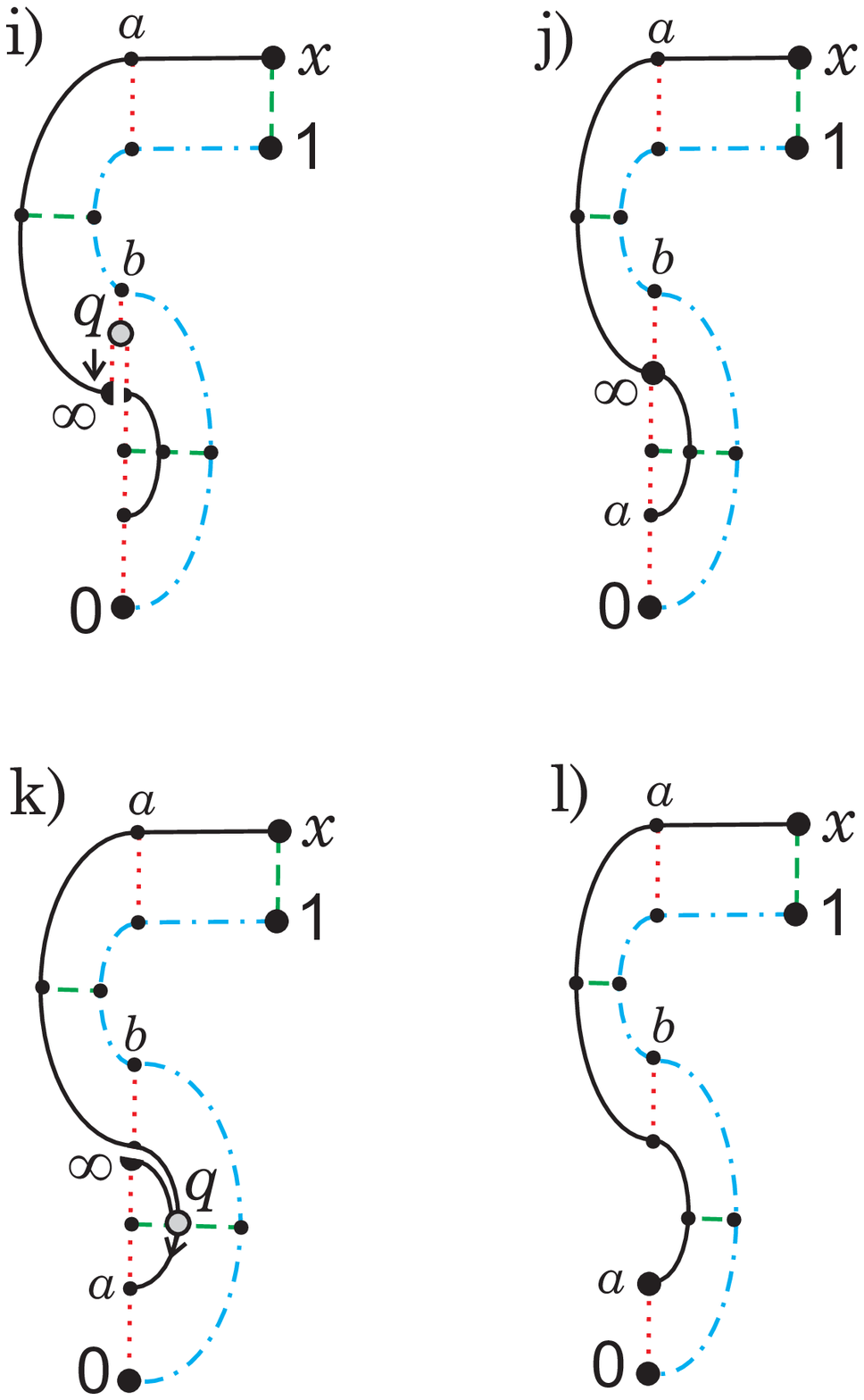}
\caption{Nets for the global family in Example 4 (continued).}\label{img12-c}
\end{figure}

\begin{figure}
\centering
\includegraphics[width=2in]{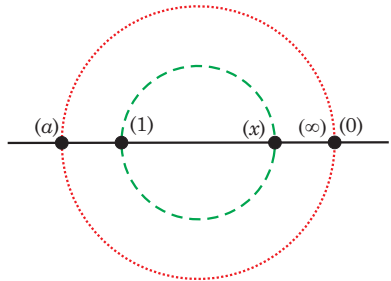}
\caption{Three-circle chain for Example 5.}\label{3circles-pent}
\end{figure}

\begin{figure}
\centering
\includegraphics[width=5in]{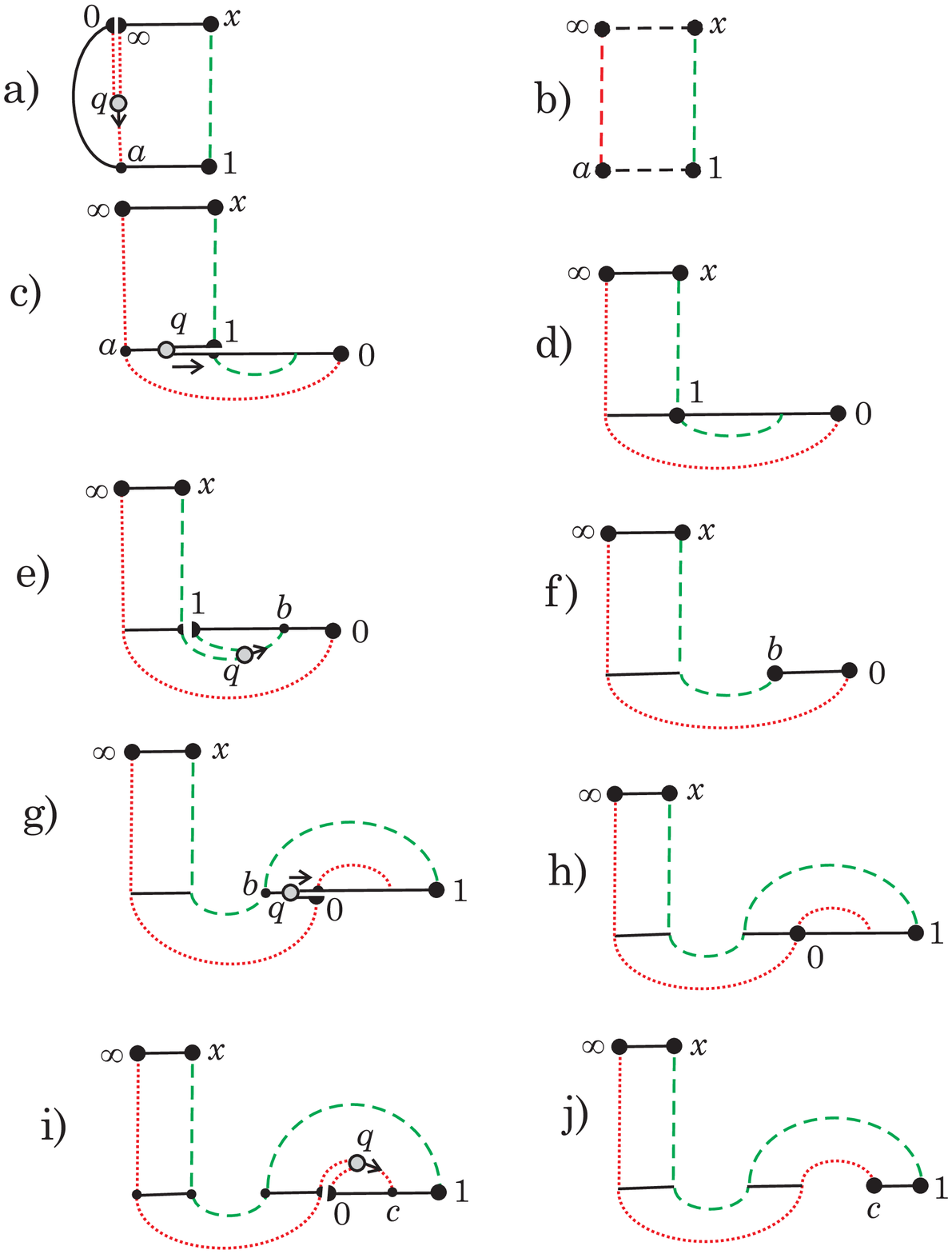}
\caption{Nets for the global family in Example 5.}\label{img11}
\end{figure}
\vspace{.1in}

\noindent{\bf Example 2.}
In this example, Fig.~\ref{pent1111}, the image of the developing map is also a region
in the sphere. In Fig.~\ref{pent1111}a, when the slit increases,
the extremal distance between $[0,1]$ and $[x,\infty]$ decreases
which means that $x$ decreases. When the slit hits the segment $[(0),(1)]$,
this distance tends to zero, which means that $x\to 1$.

As the slit in Fig.~\ref{pent1111}a decreases, $x$ increases,
and when the slit vanishes we obtain Fig.~\ref{pent1111}b.
At this moment $q(x)=x$, and we have transformation 1. As $x$ increases further
we have Fig.~\ref{pent1111}c, and then, when $(q)$ hits $[(0),(1)]$,
we obtain Fig.~\ref{pent1111}d. This is transformation 2,
and $q(x)=1$ at this point.
The digon on the right of Fig.~\ref{pent1111}c was detached.
According to the transformation 2, we attach to the quadrilateral
in Fig.~\ref{pent1111}d the vertical digon shown on
the left of Fig.~\ref{pent1111}e. The slit in Fig.~\ref{pent1111}e shortens
as $x$ increases. When it vanishes we obtain Fig.~\ref{pent1111}f
where $q(x)=0$, and we have transformation 2.
After the transformation we obtain Fig.~\ref{pent1111}g, where the slit
lengthens as $x$ increases. Eventually the slit hits $[(1),(x)]$
which corresponds to $x\to+\infty$.

Therefore, the solution $q(x)$ in this example has three special points
$x_1<x_2<x_3$ on $(1,+\infty)$,
such that $q(x_1)=x_1,\; q(x_2)=1,\; q(x_3)=0$.
These three special points correspond to quadrilaterals in
Fig.~\ref{pent1111}b,d,f. Monodromy is determined by the four circles
in Fig.~\ref{4circles-pent}b by formula (\ref{gen}).
\vspace{.1in}

\noindent
{\bf Example 3.} Consider the $4$-circle chain shown in Fig.~\ref{zeroangles1},
where all pairs $C_j,C_{j+1}$ are tangent. A quadrilateral, which is
a subset of the sphere, is the shaded region in Fig.~\ref{zeroangles2}b.
To obtain a special pentagon we make a slit $[(0),(q)]$ shown in Fig.~\ref{zeroangles2}a.
When this slit lengthens,
it eventually hits the segment $[(x),(\infty)]$ and modulus degenerates,
$x\to 1$. As the slit shortens and vanishes in Fig.~\ref{zeroangles2}b,
we have transformation 4. After that, a digon with zero
angle is attached to the shaded
region in Fig.~\ref{zeroangles2}b along a small arc $[(0),(q)]$ in
Fig.~\ref{zeroangles2}c, and the new slit $[(q),(1)]$ continues  to shorten.
The special pentagon in
this figure is not a subset of the sphere anymore: the dark shaded
area is covered twice.
When $(q)$ hits $(1)$, the slit in
Fig.~\ref{zeroangles2}c vanishes, and a new transformation 4 occurs
at a quadrilateral shown in Fig.~\ref{zeroangles2}d.
Afterwards, a sequence of transformations 4 continues indefinitely, alternately at $(0)$ and $(1)$
with the segments $[1,x]$ and $[\infty,0]$
of the pentagon increasing by a full circle
length after each two transformations,
thus $q(x)$ oscillates between $0$ and $1$ as $x\to+\infty$.
The sequence of special points is (0,1,(0,1),\ldots).
\vspace{.1in}

\noindent
{\bf Example 4.} 
Each special pentagon in this family is mapped by the developing
map to a four-circle chain shown in Fig.~\ref{4circles-pent-b}.
It is easy to check that Fig.~\ref{img12-b}b is a net corresponding
to this chain. We make a cut as shown in Fig.~\ref{img12-b}a,
and start a global family from this local family.

The nets for the global family are shown in Figs.~\ref{img12-b} and \ref{img12-c}, with the left columns (a,c,e,g,i,k)
containing the local families
of special pentagons, and the right columns (b,d,f,h,j,l)
containing quadrilaterals connecting the local families. Modulus degenerates
on one end ($x\to+\infty$). As $x\to 1$, we have an infinite chain of
local families so as $x$ increases we have the following sequence
of special points:
$$(\ldots,(\infty,\infty,0,0),\infty,x)$$
and no $1$-points. The sequence has period of length $4$, $(\infty,\infty,0,0)$
repeated infinitely many times on the left.
The monodromy in this example is
$$M_0=\left(\begin{array}{cc}0&1\\ -1&0\end{array}\right),\quad
M_1=\left(\begin{array}{cc}0&i\\ i&0\end{array}\right),\quad
M_t=\left(\begin{array}{cc}0&Ri\\ Ri&0\end{array}\right),\quad$$
where we assume that the inner and outer circles have radii $1$ and $R$.
\vspace{.1in}

\noindent
{\bf Example 5.} 
Parameters are the same as in the previous
example but monodromy is different. The global family is shown in Fig.~\ref{img11},
  with the left column (a,c,e,g,i,\dots) containing the local families
of special pentagons, and the right column (b,d,f,h,j,\dots) containing quadrilaterals connecting the local families.
The corresponding three-circle chain is shown in Fig.~\ref{3circles-pent}b.
There is an infinite
sequence of local families; as $x$ increases from $1$ to $\infty$ we have
the following sequence of special points:
$$(0,(1,1,0,0),\ldots).$$
The sequence has period $(1,1,0,0)$ repeated infinitely many times
on the right.
\begin{figure}
\centering
\includegraphics[width=2in]{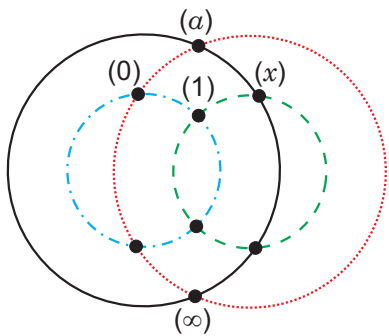}
\caption{Four-circle chain for Example 6.}\label{4circles-O}
\end{figure}
\begin{figure}
\centering
\includegraphics[width=4.5in]{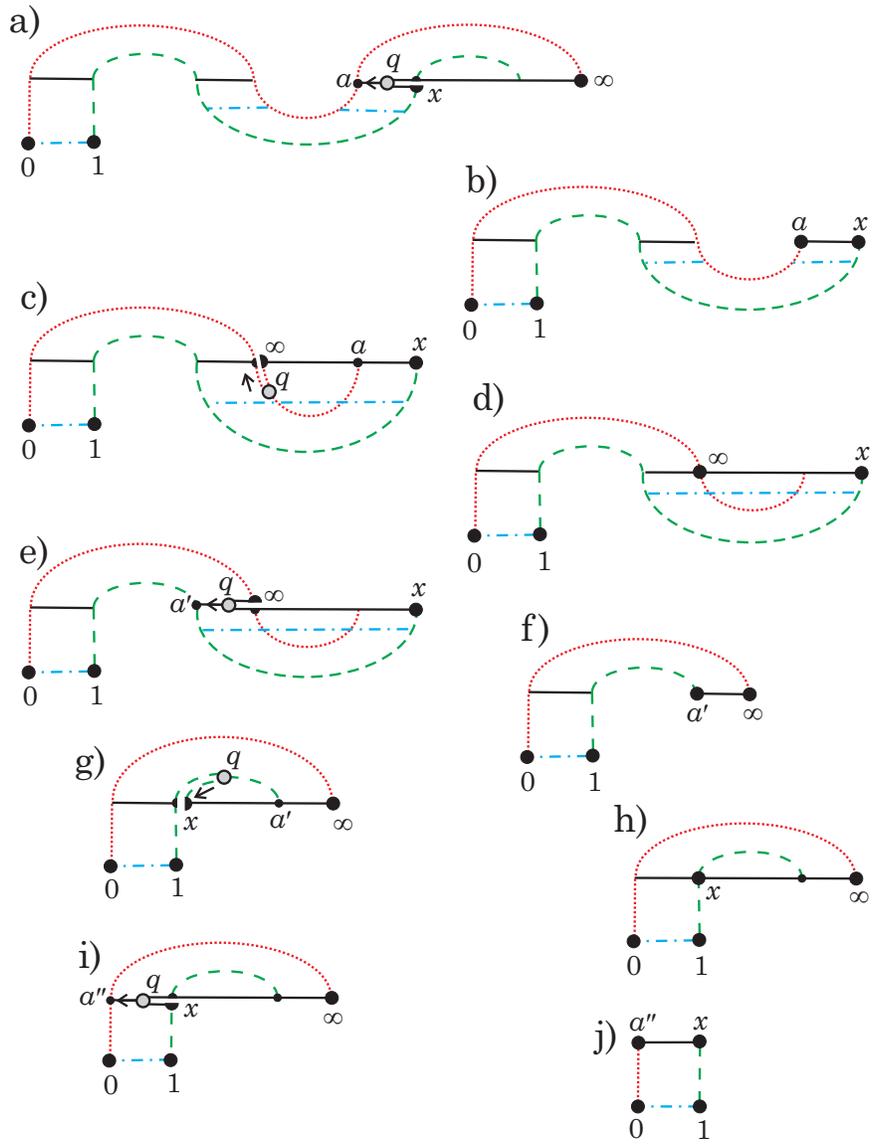}
\caption{Nets for a part of the global family in Example 6.}\label{img13a}
\end{figure}
\vspace{.1in}

\noindent
{\bf Example 6.} 
The circle chain (see Fig.~\ref{4circles-O}) consists of two pairs of non-intersecting
circles.
The global family is a doubly infinite sequence.
A part of this global family is shown in Fig.~\ref{img13a}.
It starts with a local family represented by a pentagon in Fig.~\ref{img13a}a.
As $q$ tends to $a$, this pentagon degenerates to the quadrilateral shown in Fig.~\ref{img13a}b.
In the opposite direction, when $q$ tends to $x$, the pentagon in Fig.~\ref{img13a}a does not degenerate,
and the sequence continues indefinitely, with the length of the sides $[1,x]$ and $[\infty,0]$ ever increasing.
At the other end of the sequence shown in Fig.~\ref{img13a}, a quadrilateral in Fig.~\ref{img13a}j
is symmetric with respect to reflection preserving the vertices $1$ and $a''$ and exchanging $0$ with $x$.
The sequence then continues by a local family
reflection symmetric to the pentagon shown in Fig.~\ref{img13a}i
(with the direction of $q$ reversed), and continues indefinitely,
with the length of the sides $[0,1]$ and $[x,\infty]$ ever increasing.
We have a doubly infinite sequence of special points
$$(\ldots(\infty,\infty,x,x),\infty,(0,0,\infty,\infty),\ldots).$$
The sequence is infinite in both directions,
has period $(\infty,\infty,x,x)$ on the left and $(0,0,\infty\infty)$
on the right.
\vspace{.1in}

\noindent{\bf Example 7.}
Consider the $4$-circle chain in Fig.~\ref{4circles-pent}c.
To construct a global family, we begin with a quadrilateral represented
by a net in the right column of Fig.~\ref{img15}.
That all nets in this column represent some quadrilaterals follows
from the criterion given in section 9.

Let us begin with the quadrilateral in Fig.~\ref{img15}b.
To transform it to a special pentagon, we make a slit as in
Fig.~\ref{img15}a. Lengthening of this slit corresponds to decreasing $x$.
In particular, when the slit hits the point $a$ in Fig.~\ref{img15}a,
$x\to 1$. Now we follow pictures Fig.~\ref{img15} alphabetically,
in the direction of increasing $x$. So in Fig.~\ref{img15}a the slit
shortens. As it vanishes we obtain Fig.~\ref{img15}b, transformation 1
happens, and we pass to Fig.~\ref{img15}c.

In Fig.~\ref{img15}c, the slit lengthens till $q$ hits $b$.
A transformation $2$ happens, detaching a
digon with corners $x$ and $b$ in Fig.~\ref{img15}c to obtain
the quadrilateral Fig.~\ref{img15}d. A point $c$ in Fig.~\ref{img15}d
maps to the same
point as $x$. A digon with the corners $b$ and $c$
is attached to the interval $[b,c]$ of the quadrilateral,
resulting in a pentagon Fig.~\ref{img15}e.
The slit shortens towards $x$
in Fig.~\ref{img15}e. When it hits $x$, the points $x$ and $c$ collide,
and we get a quadrilateral in Fig.~\ref{img15}f.
A transformation 1 happens at Fig.~\ref{img15}f, and the slit
lengthens towards $1$ in Fig.~\ref{img15}g. As it hits $1$, a disk with
the red (dotted line) boundary is detached, resulting in
the quadrilateral Fig.~\ref{img15}h. The point $e$ in Fig.~\ref{img15}h
maps to the same
point as $1$. A transformation 3 happens in Fig.~\ref{img15}h, with a disk
with black (solid line) boundary attached in Fig.~\ref{img15}i.
As the slit in Fig.~\ref{img15}i shortens and $q$ hits $1$, the points
$1$ and $e$ collide, and we get the quadrilateral Fig.~\ref{img15}j.

Transformations occurring in the right column of Fig.~\ref{img15}
are: 1,2,1,3,1. The last quadrilateral shown is Fig.~\ref{img15}j.
It is symmetric with respect to the reflection which exchanges
$0$ and $x$ while leaving $1$ and $\infty$ fixed.
It is easy to see that in the further continuation of the process
we will obtain all pictures Fig.~\ref{img15} in the reverse order
i)-a) subject to a reflection exchanging $0$ and $x$.

So Fig.~\ref{img15} represents only one half of the global family.
The global family is symmetric, with the symmetry exchanging $(0)$ and $(x)$.
The full sequence of special points is
$$(1,x,x,1,1,1,0,0,1).$$
\vspace{.1in}

\noindent{\bf Appendix I. Monodromy representations
corresponding to quadrilaterals}
\vspace{.1in}

\begin{figure}
\centering
\includegraphics[width=2.5in]{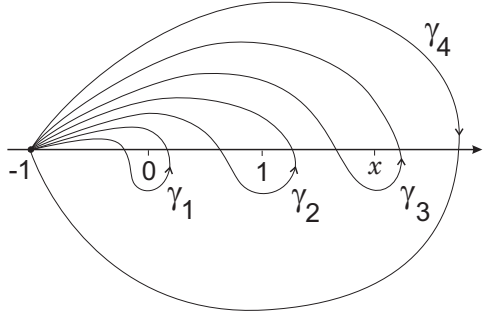}
\caption{Modified loops.}\label{loops1}
\end{figure}

\begin{figure}
\centering
\includegraphics[width=2.5in]{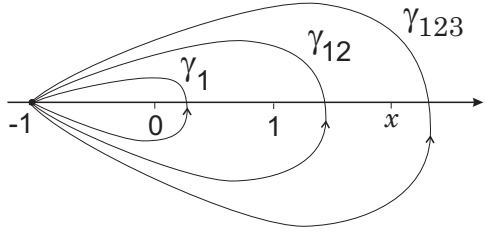}
\caption{Symmetric loops.}\label{sloops}
\end{figure}
A monodromy representation consists of $4$ matrices in $SL(2,\C)$ which satisfy
the relation (\ref{id}). For real equations (\ref{1a}) these four
matrices
can be represented as products of reflections in the circles $C_j$
containing the images of the sides of a special pentagon.
Here we will
discuss which monodromy representations correspond to real equations,
and how to find the reflections $\sigma_j$ from matrices $T_j$.

This problem was addressed in \cite{Desideri},
and we begin by restating the
result obtained there.
First of all, we change the reference point of the
fundamental group in Fig.~\ref{loops} to the point $-1$
as in Fig.~\ref{loops1}, and deform the loops accordingly.
Now consider a symmetric set of generators
of the fundamental group shown in Fig.~\ref{sloops}.
Let $N_1,N_2,N_3$ be the monodromy matrices corresponding to $\gamma_1$,
$\gamma_{12},\gamma_{123}$. We have
$$N_1=T_1,\quad N_2=T_1T_2,\quad N_3=T_1T_2T_3.$$
When none of the $\kappa_j$ is an integer, monodromy representation
determines equation (\ref{1a}) uniquely for given real $x$ and $\kappa_j$
\cite[4.2,\, 4.3]{Desideri}.
This implies that monodromy representations
correspond to real solutions of (\ref{1a})
normalized as in (\ref{normalization})
with $x_0=-1$
if and only if
\begin{equation}\label{des}
\overline{N_j}=N_j^{-1},\quad 1\leq j\leq 3,
\end{equation}
which is equivalent to the condition obtained in \cite{Desideri}.

We will derive a different condition, without the assumption
on the $\kappa_j$.

Consider an arbitrary  quadruple of $SL(2,\C)$ matrices
satisfying
\begin{equation}\label{I}
T_1T_2T_3T_4=\id.
\end{equation}
$SL(2,\C)$ acts on these quadruples by simultaneous conjugation.
To parametrize conjugacy
classes of monodromy representations we denote
$$t_j=\Tr T_j, \quad t_{jk}=\Tr(T_jT_k)=\Tr(T_kT_j).$$
Conjugacy classes are parametrized by $7$ complex numbers
\begin{equation}\label{parameters}
t_1,\; t_2,\; t_3,\; t_4,\; t_{12},\; t_{23},\; t_{13}
\end{equation}
which are subject to one relation
\begin{eqnarray}\nonumber
&t_{12}t_{23}t_{13}+t_{12}^2+t_{23}^2+t_{13}^2\\ \label{FK}
&-t_{12}(t_1t_2+t_3t_4)-t_{23}(t_2t_3+t_1t_4)-t_{13}(t_1t_3+t_2t_4) \\
&+t_1^2+t_2^2+t_3^2+t_4^2+t_1t_2t_3t_4=4.\nonumber
\end{eqnarray}
This relation was found by Fricke and Klein \cite{FK}, and was studied
in \cite{Hor}, \cite{BG} and elsewhere. Parametrization of monodromy
representations by these data is discussed in detail in \cite{iwasaki}.
In particular it is proved there that there are open dense sets on
the hypersurface (\ref{FK}) and on the space of conjugacy classes of
monodromy representations which are homeomorphic.

We say that a representation is {\em generated by reflections} if there
exist four circles $C_1,C_2,C_3,C_4$ such that the reflections
$\sigma_j$ in these circles satisfy
\begin{equation}\label{sigma}
T_j=\sigma_j\sigma_{j+1},\quad j\in\Z_4.
\end{equation}
Notice that (\ref{sigma}) implies (\ref{des}).

Arbitrary reflection can be written as
\begin{equation}\label{refl}
\sigma(z)= \frac{a\overline{z}+b}{c\overline{z}-\overline{a}},
\end{equation}
which we represent by the matrix
\def\a{\overline{a}}
\begin{equation}\label{refl2}
\left(\begin{array}{rr}a&b\\ c&-\a \end{array}\right),\quad |a|^2+bc=1,
\end{equation}
where $b,c$ are real. Product of reflections represented
by matrices $A, B$ is a linear-fractional
transformation with matrix $A\overline{B}$. Matrices $A$ associated
with reflections are characterized by the properties that $\det A=-1$ and
$\overline{A}=A^{-1}$.

Let $\Sigma_j$ be the matrices representing the reflections $\sigma_j$.
Then $\Sigma_1=I$ because of our normalization, and we have
$T_1=\overline{\Sigma_2},\; T_2=\Sigma_2\overline{\Sigma_3},\;
T_3=\Sigma_3\overline{\Sigma_4}.$
So $N_1=T_1$ has matrix $\overline{\Sigma_2},$
$N_2=T_1T_2$ has matrix $\overline{\Sigma_2}\Sigma_2\overline{\Sigma_3}
=\overline{\Sigma_3}$, and
$N_3=T_1T_2T_3$ has matrix $\overline{\Sigma_3}\Sigma_3\overline{\Sigma_4}=
\overline{\Sigma_4}.$
Thus (\ref{des}) is satisfied.

Our first question is which representations are generated by reflections.

First we notice that composition of two reflections always has real trace:
it is elliptic if the circles cross, parabolic if they are tangent and hyperbolic
if they are disjoint. Second, if (\ref{sigma}) holds then
$T_jT_{j+1}=\sigma_j\sigma_{j+2}$ also has real trace for each $j\in\Z_4.$
Thus if (\ref{sigma}) holds, the first
six parameters in (\ref{parameters}) must be real.
In addition to this we have the following inequality:
\vspace{.1in}

\noindent
{\bf Theorem A1.} {\em A monodromy representation (\ref{I}) is generated by
reflections if and only if $t_1,\; t_2,\; t_3,\; t_4,\; t_{12},\; t_{23}$
are real and
\begin{eqnarray}\label{AE}\nonumber
\Delta&:=&t_1^2t_3^2+t_2^2t_4^2+t_{12}^2t_{23}^2
-4(t_1^2+t_2^2+t_3^2+t_4^2+t_{12}^2+t_{23}^2)\\
&+& 4(t_1t_2t_{12}+t_2t_3t_{23}+t_1t_4t_{23}+t_3t_4t_{12})\\ \nonumber
&-&2(t_1t_2t_3t_4+t_2t_4t_{12}t_{23}+t_1t_3t_{12}t_{23})+16\leq 0.
\end{eqnarray}
Monodromy transformations $T_j$ determine the reflections $\sigma_j$
uniquely
unless all $T_j$ commute, and the projective monodromy group
is isomorphic to a subgroup of the multiplicative group of the unit circle
or of the additive group of the real line.
}
\vspace{.1in}

{\em Proof of Theorem A1.}

{\em Uniqueness.}  Suppose that we have (\ref{sigma}) and
\begin{equation}\label{sigmap}
T_j=\sigma_j^\prime\sigma_{j+1}^\prime,\; j\in\Z_4.
\end{equation}
First we notice that if $\sigma_j=\sigma^\prime_j$ for some $j$,
then
$\sigma_k=\sigma^\prime_k$ for all $k$. Indeed
$\sigma_j=\sigma^\prime_j$ together with (\ref{sigma}) and (\ref{sigmap})
implies $\sigma_{j+1}=\sigma^\prime_{j+1}$ and so on.

Therefore, it is sufficient to prove that $\sigma_2=\sigma^\prime_2.$
We have
\begin{equation}\label{m}
T_1=\sigma_1\sigma_2,\quad T_2=\sigma_2\sigma_3.
\end{equation}
{\bf Lemma A1.} {\em If $T_1$ and $T_2$ are non-identical
linear-fractional transformations
which together have at least three fixed points, and (\ref{m}) holds,
then $\sigma_2$ is the reflection in the unique circle which passes
through all fixed points of $T_1$ and $T_2$.}
\vspace{.1in}

{\em Proof.} If $T_1\neq\id$ then the circles $C_1$ and $C_2$
of $\sigma_1$ and $\sigma_2$
are distinct and their points of intersection are exactly the fixed points
of $T_1$. So $C_2$ contains the fixed points of $T_1$, and $T_2$.
This proves the lemma.
\vspace{.1in}

How can $T_1$ and $T_2$ have at most $2$ fixed points together?

a) One is elliptic and another one is parabolic, sharing one fixed point.

b) Both are parabolic.

c) Both are elliptic sharing two fixed points, in which case they commute.

Consider the case a). Suppose that the shared fixed point is $\infty$,
$T_1(z)=e^{2\pi i\alpha}z$ and $T_2(z)=z+c$.
Then $C_1$ and $C_2$
must be lines through the origin,
and  $C_2,C_3$  must be
parallel lines perpendicular to $c$.
Therefore  $C_2$ is the line through the origin perpendicular to
$c$, that is this circle is uniquely defined by $T_1$ and $T_2$.

Now we address b).
If in case b) $T_1$ and $T_2$ do not share their
fixed points, we may assume that
$T_1(z)=z+c$ while $T_2$ has fixed point $d\in\C$. Then $C_2$
is the unique line through $d$ perpendicular to $c$.

If the parabolic transformations in b) share the fixed point, then they
are simultaneously conjugate to $z+a$ and $z+b$,
and $C_2$ is a line perpendicular to both $a$ and $b$, so $a$ and $b$
are collinear.

So either the circle $C_2$ is uniquely defined by $T_1,T_2$,
or $T_1$ and $T_2$ commute, and either both are elliptic or both
are parabolic. If they are both parabolic, their families of invariant circles
must be the same.

This argument applies to every pair $T_k,T_{k+1}$. Therefore, the only cases
when the $\sigma_j$ are not defined by the $T_j$  are the cases stated in
the theorem. This completes the proof of uniqueness.
\vspace{.1in}

{\em Existence.}
We have already noticed
that reality of $t_1\ldots,t_4,\;t_{12},\; t_{23}$ is necessary for
(\ref{sigma}). It remains to prove that when these traces are
real, inequality (\ref{AE}) is necessary and sufficient.

We write a reflection as in (\ref{refl}), (\ref{refl2})
In particular, we obtain the
reflection in the real axis when $a=i,\; b=c=0$,
and in the line $e^{i\alpha}$ when $a=ie^{i\alpha},\; b=c=0$.

The trace of a product is
\def\Rea{\mathrm{Re}\, }
\begin{equation}\label{0}
\def\tr{\mathrm{tr}\, }
\tr(\sigma_1\sigma_2)=2\Rea(a_1\overline{a_2})+b_1c_2+c_1b_2.
\end{equation}

We normalize by $SU(2)$ conjugation so that
two adjacent circles are the real line and the line
$\{ re^{i\alpha}:r\in\R\}$, and write the four matrices
of reflections that we want to find as
$$\left(\begin{array}{rr}a_1&b_1\\ c_1&-\overline{a_1}\end{array}\right),\quad
\left(\begin{array}{rr}ie^{i\alpha}&0\\ 0&ie^{-i\alpha}\end{array}\right),\quad
\left(\begin{array}{rr}i&0\\ 0&i\end{array}\right),\quad
\left(\begin{array}{rr}a_2&b_2\\ c_2&-\overline{a_2}\end{array}\right),$$
in this order.
We can further normalize, and assume that $1$ is a fixed point
of the product of the
third and fourth reflections:
\begin{equation}\label{normalization2}
\frac{i\overline{a_2}z+ib_2}{ic_2z-ia_2}=z, \quad\mbox{where}\quad z=1,
\end{equation}
which gives
\begin{equation}\label{11}
c_2-b_2=2\Rea a_2.
\end{equation}
Now we write that the traces of products are given (real) numbers:
\begin{equation}\label{22}
2\Rea(-ie^{-i\alpha} a_1)=t_1,
\end{equation}
\begin{equation}\label{333}
2\cos\alpha=t_2,
\end{equation}
\begin{equation}\label{444}
2\Rea(i\overline{a_2})=t_3,
\end{equation}
\begin{equation}\label{55}
2\Rea(a_2\overline{a_1})+b_2c_1+c_2b_1=t_4,
\end{equation}
\begin{equation}\label{66}
2\Rea(-ia_1)=t_{12},
\end{equation}
\begin{equation}\label{77}
2\Rea(ie^{i\alpha}\overline{a_2})=t_{23}.
\end{equation}
Equations (\ref{11})--(\ref{77}) are easy to solve.
First, $a_1$ is determined from
(\ref{22}), (\ref{66}) and $a_2$ from (\ref{444}), (\ref{77}).
Then products $b_jc_j$
are found from
\begin{equation}\label{bc}
b_1c_1=1-|a_1|^2,\quad b_2c_2=1-|a_2|^2,
\end{equation}
which express the fact that determinants of our matrices are $-1$,
and together with (\ref{11}) and (\ref{55}) permits to find $b_j,c_j$.
This amounts to solving two quadratic equations. One of them always
has real solutions.
Inequality (\ref{AE}) comes from the condition that the second
also has real solutions, namely that all $b_j,c_j$ are real.

We give the details of the computation. From (\ref{22}), (\ref{66}),
\begin{equation}\label{a1}
-ia_1=\frac{1}{2}\left(t_{12}-i\frac{t_{12}\cos\alpha-t_1}{\sin\alpha}\right).
\end{equation}
Similarly, from (\ref{444}), (\ref{77}),
\begin{equation}\label{a2}
i\overline{a_2}=\frac{1}{2}\left(t_3+i\frac{t_3\cos\alpha-t_{23}}{\sin\alpha}\right).
\end{equation}
Then
\begin{equation}\label{1-a1}
1-|a_1|^2=\frac{1}{t_2^2-4}\left(t_1^2+t_2^2+t_{12}^2-t_1t_2t_{12}-4\right),
\end{equation}
and
\begin{equation}\label{1-a2}
1-|a_2|^2=\frac{1}{t_2^2-4}\left(t_2^2+t_3^2+t_{23}^2-t_2t_3t_{23}-4\right).
\end{equation}
Using (\ref{a1}) and (\ref{a2}) we obtain
\begin{equation}\label{re}
2\Rea(a_2\overline{a_1})=\frac{1}{4-t_2^2}\left(2t_1t_{23}+2t_3t_{12}-t_1t_2t_3-
t_2t_{12}t_{23}\right),
\end{equation}
and using (\ref{11})
\begin{equation}\label{c-b}
c_2-b_2=\frac{t_3\cos\alpha-t_{23}}{\sin\alpha}.
\end{equation}
Next, from (\ref{bc}), (\ref{1-a1}), (\ref{1-a2}), we obtain
\begin{equation}\label{bc1}
b_1c_1=\frac{1}{t_2^2-4}\left(t_1^2+t_2^2+t_{12}^2-t_1t_2t_{12}-4\right),
\end{equation}
and
\begin{equation}\label{bc2}
b_2c_2=\frac{1}{t^2_2-4}\left(t_2^2+t_3^2+t_{23}^2-t_2t_3t_{23}-4\right).
\end{equation}
Solving first the system (\ref{c-b}), (\ref{bc2}) with respect to $b_2,c_2$,
we obtain a quadratic equation
with discriminant
$$t_2^2t_3^2-4t_2^2-4t_3^2+16=(t_2^2-4)(t_3^2-4)\geq0,$$
because $t_j\in[-2,2].$ So we always have real solution $c_2,b_2$.

Next we solve the system (\ref{55}) with (\ref{re}) and (\ref{bc1})
with respect to $b_1,c_1$, using the known product $b_2c_2$ from (\ref{bc2}).
This also leads to a quadratic equation, whose discriminant is
a polynomial in $t_{12},\; t_{23}$ and $t_j$. This polynomial factors
(using Maple) with one factor $t_2^2-4<0$ and the other factor is $\Delta$
in (\ref{AE}).

This completes the proof.
\vspace{.1in}

{\em Remark on the proof.} In the process of recovery of $\sigma_j$
we had to solve two quadratic equations, so in general we had 4 choices to make.
On the other hand, our normalization condition (\ref{normalization2}) leaves two
choices because two circles intersect at two points.
Next, we never used $t_{13}$ in our recovery procedure for $\sigma_j$.
As $t_{13}$ satisfies the quadratic equation (\ref{FK}), assigning
$t_{13}$ narrows our choices to two.
\vspace{.1in}

An interesting question is what happens when the monodromy group is
conjugate to a subgroup of $SU(2)$.
Every element of $SU(2)$ is the product of two reflections in great circles.
If an element of $SU(2)$ is represented as a product of two reflections,
then these reflections must be in great circles, because these circles
contain the fixed points of the element which are diametrally opposite.

An interesting special case is when all seven parameters in (\ref{parameters})
are real. According to \cite[Prop. III.1.1]{MS} this happens
if and only if he projective monodromy group is a subgroup
of
$SU(2)$ or $SL(2,\R)$. When the group is generated by reflections,
the first six parameters in (\ref{parameters}) are real,
so all seven will be real if and only if the discriminant of
(\ref{FK}), as a quadratic equation with respect to $t_{13}$,
is non-negative.
Straightforward
computation shows that this discriminant is nothing but $\Delta$ defined
in (\ref{AE}). Thus we obtain
\vspace{.1in}

\noindent
{\bf Theorem A2.} {\em Let $T_1,T_2,T_3,T_4$ be unitary matrices
satisfying (\ref{I}), and all seven parameters in (\ref{parameters}
are real. This representation is generated by reflections
if and only if
$\Delta=0$, which is equivalent to
\begin{equation}\label{ass}
2t_{13}+t_{12}t_{23}-t_1t_3-t_2t_4=0.
\end{equation}
}
\vspace{.1in}

Equation (\ref{ass}) is what (\ref{FK}) becomes when $\Delta=0$.

{\em Remarks.}
We mention a simple geometric interpretation of
our conditions.

Condition that $t_{12},\; t_{23}$ are real: the trace of a product of
two elliptic transformations is real if and only if their four fixed points
lie on a circle.

Condition $\Delta=0$ gives a relation between six angles associated
to a spherical or hyperbolic quadrilateral:
four angles of the quadrilateral, and two angles between the circles
containing the images of opposite sides.
These six angles serve as natural parameters: spherical or hyperbolic
quadrilaterals
with prescribed angles at the corners form a one-parametric family, while
circular quadrilaterals with prescribed corners form a two-parametric family.

Theorem A2 can be compared with Jimbo's asymptotics \cite{Jimbo}.
(A misprint in the main result in \cite{Jimbo} was corrected in \cite{Boalch}).
It follows from the explicit formula expressing the asymptotics in
terms of the monodromy that for $SU(2)$ or $SL(2,\R)$
 monodromies this asymptotics
is real if and only if $\Delta=0$. For general monodromies (not in $SU(2)$)
it is difficult to determine directly
when Jimbo's formula gives a real asymptotics.
\vspace{.1in}

\noindent
{\bf Appendix II. Topological classification of $4$-circle chains.}
\vspace{.1in}

We recall that a $4$-circle chain consists
of $4$ labeled circles $C_j$ on the Riemann sphere such that
\begin{equation}\label{inter}
C_j\neq C_{j+1},\quad C_j\cap C_{j+1}\neq\emptyset,\quad j\in\Z_4.
\end{equation}
In this section we give
a topological classification of generic chains.
Generic means that there are no tangent circles and no triple intersections.
Two chains are
considered equivalent if there is an orientation-preserving homeomorphism
of the sphere which maps the union of circles of one chain
onto the union of circles of another chain.

\begin{figure}[H]
\centering
\includegraphics[width=5in]{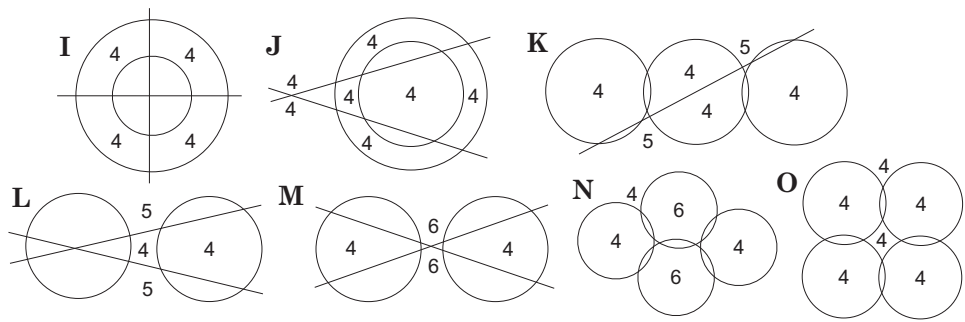}
\caption{Generic chains I-O. Some pairs are disjoint.}\label{chains2}
\end{figure}

In fact, we classify generic {\em unordered} quadruples of circles with the
following property: each circle intersects at least two other circles.
There are two kinds of such quadruples: those in which each circle intersects all three
other circles (see Fig.~\ref{chains1}) and those with a pair of non-intersecting circles
(see Fig.~\ref{chains2}). Note that the quadruples in Fig.~\ref{chains1}D and \ref{chains2}K
are not reflection symmetric, thus each of them represents two equivalence classes.
We omit the elementary but tedious proof that these exhaust all possibilities.
To see that all these cell decompositions are
distinct we indicated the faces with more than 3 edges in each
cell decomposition.

It follows from the classification that the circles in Fig.~\ref{chains1} can be
arbitrarily ordered to form a $4$-circle chain, while the circles in Fig.~\ref{chains2}
form a $4$-circle chain when ordered so that non-intersecting circles are not adjacent.

Notice two equivalent reformulations of this problem:
topological classification of arrangements of
four planes in the hyperbolic space, subject to the intersection
condition (\ref{inter}), and topological classification of
possible intersections of a sphere with four planes
in the Euclidean space, under the condition that the planes $P_j$
can be so ordered that each line $P_j\cap P_{j+1}$ intersects
the sphere.

\vspace{.1in}

{\em Remarks and conjectures.}
\vspace{.1in}

If all pairs of circles in the chain intersect, then there are only finitely
many nets on this chain with prescribed angles. This follows from the results
of \cite{Ihlenburg}. In this paper, Ihlenburg proves that all circular
quadrilaterals can be obtained from finitely many topological types by
applying four explicitly defined operations. All these operations
do not decrease angles, and three of them
increase some angles. The only operation which leaves all angles unchanged
requires two disjoint circles in the chain.

It follows that real solutions of PVI corresponding to all chains
in Fig.~\ref{chains1} can have only
finitely many special points on an interval between fixed singularities.
On the other hand, our examples 4, 6 suggest that for all
chains containing pairs of disjoint circles
the number of special points is infinite.
Moreover, it looks like it is infinite in one direction when there is one
pair of disjoint circles, and infinite in both directions
if there are two such pairs, like in Fig.\ref{4circles-O} which is the same
as Fig.~\ref{chains2}O.

\nopagebreak

\noindent
{\bf Appendix III. Moduli of conformal quadrilaterals.}
\vspace{.1in}

Consider a closed rectangle $Q$ in the plane with vertices $0,1,1+ia,ia$.
\def\mod{\mathrm{mod}\ }
The number $a>0$ is called the modulus, $a=\mod Q$.
Any Borel measurable function $\rho(z)\geq 0$ defines a {\em conformal metric}
$\rho(z)|dz|$
on $Q$: the length of a curve $\gamma$ and the area of a set $E\subset Q$
are defined as
$$\ell_\rho(\gamma)=\int_\gamma\rho(z)|dz|, \quad\mbox{and}\quad
A_\rho(E)=\int_E\rho^2(z)dxdy.$$
Let $\Gamma$ be the set of all curves in $Q$ connecting the horizontal
sides. Define
$$\ell_\rho(\Gamma)=\inf_{\gamma\in\Gamma}\ell_\rho(\gamma),$$
and
\begin{equation}\label{111}
\lambda(\Gamma)=\sup_\rho\frac{\ell_\rho^2(\Gamma)}{A_\rho(Q)},
\end{equation}
where the sup is taken over all metrics for which the numerator and denominator
are finite and not zero.
\vspace{.1in}

\noindent
{\bf Lemma A2.} \cite[I.D, Example 1]{Ahl2} $\lambda(\Gamma)=a.$
\vspace{.1in}

Formula (\ref{111}) defines the {\em extremal length} of
an {\em arbitrary} family of curves $\Gamma$
in $Q$. So defined extremal length is a conformal invariant
of a family of curves.

Let $\Gamma'$ be the family of all curves in $Q$ connecting
the vertical sides. Then evidently
\begin{equation}\label{2}
\lambda(\Gamma)\lambda(\Gamma')=1.
\end{equation}
The following comparison inequalities immediately follow the from
definition.
\vspace{.1in}

\noindent
{\bf Lemma A3.} \cite[I.D, Theorem 2]{Ahl2} {\em Consider two families
of curves $\Gamma_1$ and $\Gamma_2$ and suppose that every curve $\gamma_2\in \Gamma_2$ contains some curve $\gamma_1\in\Gamma_1$. Then
$\lambda(\Gamma_1)\leq\lambda(\Gamma_2)$.}
\vspace{.1in}

The assumption means that $\Gamma_1$ has ``more curves'' and the curves of
$\Gamma_2$ are ``longer''.

For a metric $\rho$, the intrinsic distance $d_\rho(E_1,E_2)$
between two subsets $E_1$ and $E_2$
of $Q$ is defined as infimum of $\ell_\rho(\gamma)$ over all curves
connecting a point in $E_1$ with a point in $E_2$.
\vspace{.1in}

\noindent
{\bf Lemma A4.} {\em
Suppose that a metric $\rho$ has the following properties:
\begin{equation}\label{area}
A_\rho(B(r))\leq Kr^2,
\end{equation}
for all $r>0$ and for all intrinsic disks $B(r)$ of radii $r$,
the intrinsic $\rho$-distance between the vertical sides is at least $2c$,
and
the intrinsic $\rho$-distance between the
horizontal sides is less than $\epsilon/2$.

Then $\mod Q\leq \delta$, where
$$\delta=\frac{4(K+1)}{\log(c/\epsilon)}\to 0$$
as
$\epsilon\to 0$, for any fixed $K>0,c>0$.}
\vspace{.1in}

{\em Proof.} Choose a curve $\gamma_0$ connecting the horizontal sides,
and such that $\ell_\rho(\gamma_0)<\epsilon$. Let $P$ be a point on $\gamma_0$.
Consider the closed $\rho$-disks
$B(\epsilon)$ and $B(c)$ of radii $\epsilon$ and $c$,
both centered at $P$. Then $B(\epsilon)$ contains $\gamma_0$.
Let $\Gamma'$ be the family of all curves in
$Q$ connecting the vertical sides.
Every curve $\gamma'$ of this family crosses $\gamma_0$, therefore
$\gamma'$ intersects both $B(\epsilon)$ and $Q\backslash B(c)$.
Therefore
\begin{equation}\label{3}
\lambda(\Gamma')\geq\lambda(\Gamma_1),
\end{equation}
where $\Gamma_1$ is the family of all curves in $Q$ connecting
$B(\epsilon)$ with $Q\backslash B(c)$. To estimate $\lambda(\Gamma_1)$
from below, consider the metric $\tau(z)|dz|$ defined by the function
$$\tau(z)=\frac{\rho(z)}{\log(c/\epsilon) d_\rho(z,P)},
\quad
z\in B(c)\backslash B(\epsilon),$$ and zero otherwise.
For every $\gamma_1\in \Gamma_1$ we have
$$\ell_\tau(\gamma_1)\geq
\frac{1}{\log(c/\epsilon)}\int_{\epsilon}^{c}
\frac{\rho(z)|dz|}{d_\rho(z,P)}
\geq
\frac{1}{\log(c/\epsilon)}\int_{\epsilon}^{c}\frac{ds}{s}\geq
1.$$
Here we made the change of the variable $s=d_\rho(z,P)$
and used the evident inequality $ds\leq\rho(z)|dz|$.

To estimate the $\tau$-area of $Q$ we define $r_0=\epsilon,\; r_k=2^kr_0$,
$k=1,\ldots,N,$ $N=[\log_2(c/\epsilon)]+1$, and let $B_k$ be the $\rho$-disk
of radius $r_k$ centered at $P$.
Then, using (\ref{area}), we obtain
\begin{eqnarray*}
A_\rho(Q)&\leq & \frac{1}{\log^2(c/\epsilon)}\sum_{k=1}^N
\int_{B_k\backslash B_{k-1}}\frac{\rho^2(z)dxdy}{r_{k-1}^2}\leq
\frac{4}{\log^2(c/\epsilon)}\sum_{k=1}^N\frac{A_\rho(B_k)}{r_k^2}\\
&\leq &
\frac{4KN}{\log^2(c/\epsilon)}\leq \frac{4(K+1)}{\log(c/\epsilon)}.
\end{eqnarray*}
Thus $\lambda(\Gamma_1)\geq (\log(c/\epsilon))/(4(K+1))$,
and using (\ref{2}) and  (\ref{3}), we obtain
$$\mod Q=\lambda(\Gamma)\leq 4(K+1)/(\log(c/\epsilon)).$$
This proves the lemma.
\vspace{.1in}

The upper half-plane can be mapped conformally onto $Q$ so that
$$(0,1,x,\infty)\mapsto (0,1,1+ia,ia)$$
by the Schwarz--Christoffel formula. Let
$$\phi(z)=\int_0^z\frac{d\zeta}{\sqrt{z(z-1)(z-x)}}.$$
Then the desired conformal map is $\phi(z)/\phi(1)$,
and the modulus $a(x)=-i\phi(\infty)/\phi(1).$
It follows from Lemma~A3 that $x\mapsto a(x)$ is increasing homeomorphism
of $(1,\infty)$ onto $(0,+\infty)$.
\vspace{.1in}

In our applications, the metric $\rho$ arises as a pull-back of
the standard spherical metric of curvature $1$ on the sphere $S$
by a conformal local homeomorphism $f:Q\to S$.
If $f$ is $p$-valent (which means that every point has at most $p$
preimages), then
(\ref{area}) is satisfied with $K=\pi p$. Indeed, for spherical discs
on $S$, (\ref{area}) is satisfied with $K=\pi$ by direct computation,
and $f(B(r))$ is evidently contained in a disk of radius $r$
in $S$.

\vspace{.1in}

{\em Department of Mathematics, Purdue University,
West Lafayette, IN 47907 USA}
\end{document}